\documentclass[11pt]{article}

\usepackage{amsmath,amsfonts,amssymb}
\usepackage{graphicx}
\usepackage{color}
\usepackage{caption}

\def\R{{\mathbb R}}

\def\L{{\mathcal L}}

\let\theta\vartheta
\let\phi\varphi

\def\vec#1{{\mathbf{#1}}}

\def\uh{\vec{u}_h}
\def\eh{\vec{e}_h}
\def\fh{\vec{f}_h}

\def\rh{\vec{r}_h}
\def\wh{\vec{w}_h}
\def\uhh{\vec{u}_{2h}}
\def\ehh{\vec{e}_{2h}}

\def\rhh{\vec{r}_{2h}}

\newcommand{\pad}[2]{\frac{\partial #1}{\partial #2}}

\title{Second Order Multigrid Methods for Elliptic Problems with Discontinuous Coefficients on an Arbitrary Interface, I: One Dimensional Problems}
\author{Armando Coco, Giovanni Russo\thanks{Dipartimento di Matematica e Informatica, Universit\`a di Catania, Catania, Italy}}
\date{}

\begin{document}

\maketitle

\begin{abstract}
In this paper we present a one dimensional second order accurate method to solve Elliptic equations with discontinuous coefficients on an arbitrary interface. Second order accuracy for the first derivative is obtained as well. The method is based on the Ghost Fluid Method, making use of ghost points on which the value is defined by suitable interface conditions. The multi-domain formulation is adopted, where the problem is split in two sub-problems and interface conditions will be enforced to close the problem. Interface conditions are relaxed together with the internal equations (following the approach proposed in~\cite{CocoRusso:Elliptic} in the case of smooth coefficients), leading to an iterative method on all the set of grid values (inside points and ghost points). A multigrid approach with a suitable definition of the restriction operator is provided. The restriction of the defect is performed separately for both sub-problems, providing a convergence factor close to the one measured in the case of smooth coefficient and independent on the magnitude of the jump in the coefficient.
Numerical tests will confirm the second order accuracy.
\\Although the method is proposed in one dimension, the extension in higher dimension is currently underway~\cite{CocoRusso:discontinuous2d} and it will be carried out by combining the discretization of~\cite{CocoRusso:Elliptic}
with the multigrid approach of~\cite{CocoRusso:MG} for Elliptic problems with non-eliminated boundary conditions in arbitrary domain.
\end{abstract}

\section*{Introduction}
Elliptic equations with jumping coefficients across a one-codimensional interface $\Gamma$ arise in several applications. Let us mention as examples the steady-state diffusion problem in two materials with different diffusion coefficient separated by an arbitrary interface, the Poisson equation coming from the projection method in incompressible Navier-Stokes equation for fluids with different density, the porous-media equation to model the oil reservoir, electrostatic problems, and many others. In order to close the problem, interface conditions related to the jump of the solution and of the flux across the interface are included.
In all these problems the interface may be arbitrary (not aligned with a line grid) and can change in time.

Numerous techniques have been developed to treat such problem. Interface-fitted grid methods such as the ones based on Finite Elements Methods~\cite{Babuska:FEM, Bramble:FEM}
are not suitable in case of moving interface, because a re-meshing grid is needed at each time step and this makes the computation expensive. Then an approach treating the interface embedded in a Cartesian grid and moving according to the velocity field of the fluid is preferred. Since the interface may not be aligned with the grid, a special treatment is needed. The simplest method makes use of the Shortley-Weller discretization~\cite{Shortley-Weller:discretization}, that discretizes the Laplacian operator with usual central difference away from the interface, and makes use of a non symmetric stencil in the points close to the interface, adding extra-grid points on $\Gamma$. While jumping condition on the solution is straightforward to discretize on interface points, the jump in the flux (involving the normal derivative) cannot be immediatly discretized in more than one dimension. In fact, Shortley-Weller discretization requires that the value of the normal derivative of the solution on both sides of the interface is suitably reconstructed at the intersection between the grid and the interface. This approach is adopted, for example, by Hackbusch in~\cite{Hackbusch:elliptic} to first order accuracy, and by other authors (see~\cite{Bramble:ell} and references therein) to second order accuracy. However, the method proposed by Bramble in~\cite{Bramble:ell} for second order accuracy is quite involved and not recommendable for practical purposes.

Methods based on embedding the domain in a Cartesian grid without adding extra-grid points are derived from the pioneering work of Peskin~\cite{Peskin:IBM}, where the Immersed Boundary Methods is introduced to model blood flows in the heart. In that paper a source term is localized on the the boundary and the method makes use of a discretized delta-function, leading to a first order accuracy. 
A second order accurate extension to jump coefficients is the Immersed Interface Methods, first developed by LeVeque and Li in~\cite{LeVequeLi:IIM}. Such method uses a six-point stencil to discretize the elliptic equation in grid points close to the interface $\Gamma$ and the coefficients of such stencil are found by Taylor expansion of the solution. Jump conditions on the interface are then used to modify the coefficients appearing in the equation corresponding to nodes near $\Gamma$, in such a way that the overall discretization is second order accurate.
Non-homogeneous jump conditions are allowed on the function and on the normal flux. 

Another method which achieves second order accuracy by modifying standard difference formulas was proposed by Mayo in~\cite{Mayo:biharmonic} for solving Poisson or biharmonic equation on irregular domains. Such method embeds the irregular domain in a regular region with a Cartesian grid and discretizes the equation on the whole region, by suitable extension of the solution outside.

In all these methods the only unknowns are the values on the grid points and the stencil may cross the interface, leading to a quite involved procedure to reach the desired accuracy, since the derivative of the solution may jump crossing the interface and values from the other side are used in the computation.

A rather simple method to use standard five-point stencil even close to the interface is the Ghost-Fluid Method, introduced by Fedkiw \textit{et al.} in~\cite{Fedkiw:GFM}. Here the authors point out that a two-phase problem could be reduced in two sub-problems by a multi-domain formulation, and each sub-problem may be discretized with the same technique used to solve a single problem with Dirichlet/Neumann boundary conditions. Such method makes use of extra grid points (\textit{ghost points}) outside the domain in order to keep unchanged the symmetry of the stencil even for inside points close to the interface. In ghost points, interface conditions are enforced in order to close the discrete system.

Methods based on ghost points are discussed in~\cite{Gibou:Ghost}, where Gibou \emph{et al.} proposed a second-order accurate method for Dirichlet conditions on regular Cartesian grid. The value at the ghost nodes is assigned by linear extrapolation, and the whole discretization leads to a symmetric linear system, easily solved by a preconditioned conjugate gradient method. A fourth order accurate method is also proposed in~\cite{Gibou:fourth_order}. Other methods use a non-regular Cartesian grid, such as in~\cite{Gibou:quadtree}, where Gibou \emph{et al.} present finite difference schemes for solving the variable coefficient Poisson equation and heat equation on irregular domains with Dirichlet boundary conditions, using adaptive Cartesian grids. One efficient discretization based on cut-cell method to solve more general Robin conditions is proposed by Gibou \textit{et al.} in~\cite{Gibou:Robin}, which provides second order accuracy for the Poisson and heat equation and first order accuracy for Stefan type problems.

Other approaches based on cut-cell methods obtained by a Finite Volume discretization are presented in~\cite{Colella:PoissonFV}. Cells that are cut by the boundary requires a special treatment, such as cell-merging and rotated-cell, in order to avoid a too  strict restriction of the time step dictated by the CFL condition.

Several methods have been also proposed to model the interaction between multiphase flows and solid obstacles, such as Arbitrary Lagrangian Eulerian (ALE)~\cite{FormaggiaNobile:ALE, Donea:ALE}, Distribute Lagrangian Multiplier (DLM)~\cite{Glowinski:DLM}, penalization methods~\cite{Sartou:penalization, Angot:penalization}. In~\cite{Iollo:penalization} a combination of penalization and level-set methods is presented to solve inverse or shape optimization problems on uniform Cartesian meshes.
In~\cite{Zhou:MIB} Zhou \textit{et al.}  proposed a Matched Interface and Boundary (MIB) method for elliptic problems with sharp-edged interfaces. 

In time-dependent problems requiring the solution of an elliptic problem at each time step an iterative solver is
preferred with respect to a direct problem, since a good initial guess (the solution at the previous time step) is
provided. Most iterative method for jumping coefficient are based on Domain Decomposition
Methods~\cite{QuarteroniValli:DDM}, either with or without overlapping. Such methods are based on the multi-domain
formulation, i.e., the problem is split in two sub-problems and interface conditions are enforced to achieve two
sub-problems with respectively Dirichlet and Neumann boundary (coupled) conditions on the interface/boundary. Each
sub-problem is solved and the solution at the interface is used to provide an updated right-hand side for the other
sub-problem, and so on iteratively. A drawback of this method is that association between the Dirichlet/Neumann boundary
condition and the sub-domain cannot be arbitrary (see~\cite[pag.~12]{QuarteroniValli:DDM}).

Most applications require second order accuracy in the gradient: for example, in projection method for incompressible Navier-Stokes equation, the gradient of the pressure is used to correct the fictitious velocity field leading it to satisfy the free-divergence condition. Also high-order accuracy~\cite{Gibou:fourth_order} may be required, for instance when turbulence and shock interact, or high frequency wave propagation are presented in inhomogeneous media~\cite{Zhao:HighOrder}.

In~\cite{CocoRusso:Elliptic} a second-order accurate discretization for elliptic problems in arbitrary domain and mixed boundary condition is provided, together with a convergence proof for the iterative solver for first order accuracy. The method is based on transforming the stationary problem into a fictitious evolutionary problem, both inside the domain and on the boundary. The problem is then discretized on a regular grid using non eliminated boundary conditions to determine the proper relaxation equation for the ghost points. The whole procedure is made efficient by a multigrid technique, as illustrated in 
\cite{CocoRusso:MG}.

The present paper provides a second order discretization of the problem based on the ghost-point method on regular Cartesian grid described in~\cite{CocoRusso:Elliptic} and makes use of an iterative solver whose convergence is speeded up by a multigrid approach~\cite{CocoRusso:MG}. Interface conditions are neither eliminated from the discrete system (they are strongly coupled and their elimination is too hard to perform in more than one dimension) nor directly enforced (which leads to a non-convergent iterative method): they are relaxed together with the interior equations. This leads us to an iterative scheme for the set of all unknowns (internal points and ghost points). The method works also for non-homogeneous interface conditions.
Although this paper provides a 1D description of the method, the generalization of the approach in higher dimension is currently underway~\cite{CocoRusso:discontinuous2d} and can be obtained in an almost straightforward manner combining results from~\cite{CocoRusso:Elliptic, CocoRusso:MG}.

Several multigrid approaches exist in literature to treat the jumping coefficient problem in 2D when the interface is aligned with the Cartesian grid. We mention the method based on operator-dependent interpolation~\cite{Brandt:DiscCoeff, Dendy:BlackBoxMG}, where the interpolation is carried out by exploiting the continuity of the flux instead of the gradient of the solution, and the method based on Galerkin Coarse Grid Operator~\cite{Stuben:AMG}, which makes the algebraic problem more expensive from a computational point of view and does not take advantage from the fact that the discrete problem comes from a continuous problem.

In our approach we use the standard interpolation operator and discretize the operator in the coarser grid in the same way as in the fine grid, without making use of Galerkin conditions.
But, since the defect may jump crossing the interface, a separated restriction for both sub-problems is needed, as performed in~\cite{CocoRusso:MG} for arbitrary domain with mixed boundary condition (without jumping coefficient). This approach provides a good convergence factor, comparable with ones measured for no-jumping case. We also show that the convergence factor does not depend on the magnitude of the jump in the coefficient. Interface conditions are relaxed, then have to be transferred to the coarse grid as well. In one-dimensional case this task is trivial, since such conditions are just two real values that can be copied to the coarse grid. In higher dimension interface conditions are stored in ghost points, which can show a complex structure for arbitrary interface. The restriction of interface condition defect can be carried out in the same manner of the restriction of boundary condition defect described in~\cite{CocoRusso:MG} for problems with non-eliminated boundary conditions: the defect is first extrapolated outside the domain and then transferred to the coarse grid in the same manner as the restriction of the defect of inside equations, i.e., without using values from the other side of the boundary. This work is currently underway.

The rest of the paper is divided in 3 sections. In the first section we describe the second order accurate discretization of the model problem and the iterative scheme obtained by the relaxation of the interface conditions.
The second section is devoted to the multigrid approach, with a care description of the transfer operators. In section 3 some numerical test is performed, to show the second order accuracy in the solution and in its first derivative as well. We measure also the convergence factor and compare it with the convergence factor obtained by other methods. 

\section{Second order accurate discretization}
In this section we obtain a second order accurate numerical method to solve an elliptic equation with discontinuous coefficients. After introducing the model problem, we provide a discretization and an iterative solver of the linear system. In some applications one may be interested in second order accuracy also for the derivative of the solution. In numerical tests of Sec. \ref{NumTests} we show that the method is second order accurate in the solution and in its first derivative.

\subsection{Model problem} 	
Let us consider the model problem
\begin{equation}\label{modprob}
 \begin{split}
-   \displaystyle \frac{d}{dx} \left( \gamma \displaystyle \frac{d u}{dx} \right) =& f \mbox{ in } \Omega=[0,1],\\
u(0)=g_0, \; \; \; \; \; \; & u(1)=g_1.
\end{split}
\end{equation}
where the diffusion coefficient $\gamma \colon [0,1] \rightarrow \R$ jumps on an interface $\alpha \in ]0,1[$, i.e., is a smooth function in $[0, \alpha[$ and in $]\alpha, 1]$, but may be discontinuous across $\alpha$. We assume $\gamma>\epsilon>0$ in all the domain.
If we solve this problem by standard central differences on a uniform grid, the accuracy of the method degrades to first order.

\begin{figure}[!hbt]
 \centering
   	\includegraphics[width=0.60\textwidth]{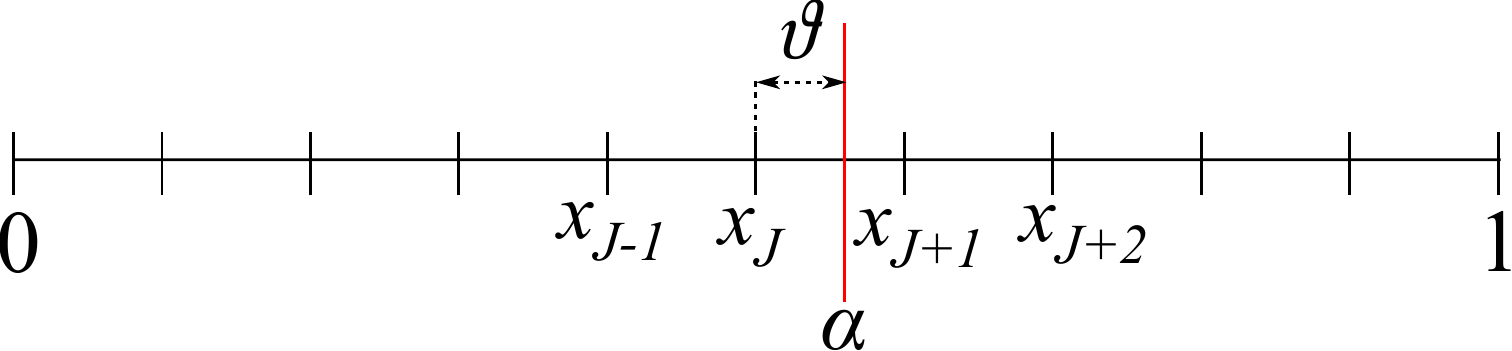}
   	\captionsetup{width=0.80\textwidth}
	\caption{\footnotesize{ Computational domain $\Omega$ with an arbitrary interface $\alpha$. }}
	\label{fig:Omega1d}
 	\end{figure}
 	
Let
\[
u^L = u|_{[0,\alpha[}, \; \; \; u^R = u|_{]\alpha,1]}, \; \; \; \gamma^L = \gamma|_{[0,\alpha[}, \; \; \; \gamma^R = \gamma|_{]\alpha,1]}
\]
be the restriction functions of the solution and of the coefficient on the two subdomains.
We split the problem into the following subproblems:
\begin{equation}\label{sub1}
 \begin{split}
- \displaystyle \frac{d}{dx} \left( \gamma^L \displaystyle \frac{du^L}{dx} \right) &= f \mbox{ in } [0,\alpha[ \\
u^L(0)&=g_0,
\end{split}
\end{equation}
\begin{equation}\label{sub2} 
\begin{split}
- \displaystyle \frac{d}{dx} \left( \gamma^R \displaystyle \frac{d u^R}{dx} \right) &= f \mbox{ in } ]\alpha,1] \\
u^R(1)&=g_1.
\end{split}
\end{equation}
In order to close the problem, we must provide an additional boundary condition for each of $u^L$ and $u^R$ on the interface $\alpha$. This additional conditions are inferred to the requirement that the solution $u$ and the flux $\gamma u'$ are continuous across $\alpha$. Introducing the \textit{jumping operator} on $\alpha$
\[
\left[ w \right] = \lim_{x \rightarrow \alpha^+} w - \lim_{x \rightarrow \alpha^-} w,
\]
the additional boundary conditions may be resumed as
\[
\left[ u \right] = 0, \; \; \; \left[ \gamma u' \right] = 0
\]
and are called \textit{transmission conditions}~\cite{QuarteroniValli:DDM}. They can be inferred by a physical requirement: for instance, in steady-state diffusion problems in two materials, the temperature and its flux are required to be continuous across $\alpha$. Non-homogeneous interface conditions may appear, for example, in presence of a delta-function on the right hand side $f=f_1+\delta_\alpha$, with $f_1 \in C^0([0,1])$. Precisely, the two following problems are equivalent:
\begin{equation*}
\begin{split}
- \displaystyle \frac{d}{dx}(\gamma \, \displaystyle \frac{du}{dx})=&f_1+C\delta_\alpha \mbox{ in } [0,1] \\
u(0)=g_0, \; \; \; \; \; \; & u(1)=g_1,
\end{split}
\end{equation*}
\begin{equation*}
\begin{split}
-\displaystyle \frac{d}{dx}(\gamma \: \displaystyle \frac{du}{dx})=&f_1 \mbox{ in } [0,\alpha[ \cup ]\alpha,1]\\
u(0)=g_0, \; \; \; \; \; \; & u(1)=g_1 \\
\left[ u \right]=0, \; \; \; \; \; \; &  \left[ \gamma \: u' \right]=-C.
\end{split}
\end{equation*}
In the following we suppose the right-hand side is a regular function in the two sub-regions, and non-homogeneous interface conditions are allowed:
\begin{equation}\label{transmission}
\left[ u \right] = g_D, \; \; \; \left[ \gamma u' \right] = g_N.
\end{equation}
Such general case is relevant for some applications, for example pressure equation for incompressible flow in presence of surface tension at the interface.
\\The two subproblems (\ref{sub1}) and (\ref{sub2}) are then coupled on $\alpha$ and cannot be solved separately. The whole problem becomes
\begin{eqnarray}
- \displaystyle \frac{d}{dx} \left( \gamma^L \displaystyle \frac{d\,u^L}{dx} \right) &=& f \mbox{ in } [0,\alpha[ \label{eqL} \\
- \displaystyle \frac{d}{dx} \left( \gamma^R \displaystyle \frac{d\,u^R}{dx} \right) &=& f \mbox{ in } ]\alpha,1] \label{eqR} \\
u^L(0)=g_0, && u^R(1)=g_1 \label{DirBC}\\
\left[ u \right] = g_D, && \left[ \gamma u' \right] = g_N. \label{jumpsU_Flux}
\end{eqnarray}

\subsection{Discretization}\label{discr}
Let $N$ be an integer, $h=1/(N+1)$ be the spatial step and $x_0, x_1, \ldots, x_{N}, x_{N+1}$ be the equally spaced grid points, with $x_j=j\,h$. Let $J$ be such that $x_J \leq \alpha < x_{J+1}$ (see Fig. \ref{fig:Omega1d}). We write $J=\left\lfloor \alpha \right\lfloor$, where $\left\lfloor \cdot \right\lfloor$ denotes the integer part. We will denote by $\L_{j}[w]$ the quadratic interpolant of $w$ in nodes $\left\{ x_{j-1}, x_j, x_{j+1} \right\}$. 
By $u_j^L$ [$u_j^R$] we denote the component of the numerical solution which approximates $u^L(x_j)$ [$u^R(x_j)$], while we intend $f_j=f(x_j)$, $\gamma^L_j=\gamma^L(x_j)$, $\gamma^R_j=\gamma^R(x_j)$. 
\\Let us discretize the system (\ref{DirBC}). 
Discretizing Eq. (\ref{eqL}) on nodes $x_1, x_2, \ldots, x_{J} $ using central differences for the solution $u^L$ and linear interpolation for the coefficient function $\gamma^L$, we obtain:
\begin{equation}\label{eqLj}
\displaystyle \frac{1}{h^2} \left( \gamma^L_{j- \frac{1}{2}} \left( u^L_j-u^L_{j-1} \right) + \gamma^L_{j+ \frac{1}{2}} \left( u^L_j-u^L_{j+1} \right) \right) = f_j, \; \; j=1, \ldots J,
\end{equation}
where $\gamma^L_{j + 1/2}=(\gamma^L_{j}+\gamma^L_{j + 1})/2$. In Eq. (\ref{eqLj}) for $j=1$ the value $u^L_0$ is given by the Dirichlet condition (\ref{DirBC}): $u^L_0=g_0$.
It can be easily eliminated from (\ref{eqLj}), but we will leave it in the system just for simplicity. The same applies for $u^R_{N}$ discretizing Eq. (\ref{eqR}) in node $x_{N}$.
\\Eq. (\ref{eqLj}) for $j=J$ needs to know the values of $u^L$ and $\gamma^L$ in node $x_{J+1}$. Since $u'$ and $\gamma$ are discontinuous, we cannot use respectively $u^R_{J+1}$ and $\gamma^R_{J+1}$, because this may result in a loss of accuracy, since it smears out the coefficient $\gamma$ and the numerical solution itself, while both jump on the interface. Then we need to add an additional grid point value for the numerical solution $u^L(x_{J+1})$, called \textit{ghost point value}, and to extrapolate $\gamma^L$ up to the first ghost point $x_{J+1}$. The same argument holds for $u^R$ and $\gamma^R$ in their ghost point $x_J$, when discretizing Eq. (\ref{eqR}) in node $x_{J+1}$.
\\The unknowns of the numerical method are therefore the $N+4$ quantities
\begin{equation}\label{unknowns}
u_0^L, \ldots, u_{J+1}^L, u_{J}^R, \ldots, u_{N+1}^R.
\end{equation}
This approach has been called \textit{Ghost Fluid Method} and used in the context of multi-fluid flows~\cite{Fedkiw:GFM}. 
The two additional unknowns $u_{J+1}^L$ and $u_{J}^R$ require two additional boundary conditions to close the system, which are given by the transmission conditions (\ref{transmission}), resulting in a $2 \times 2$ sub-system. We will not solve this sub-system for $u_{J+1}^L$ and $u_{J}^R$, but we instead leave it in the whole linear system, which will be solved iteratively.
The extrapolation for the coefficient functions $\gamma^L$ and $\gamma^R$ is simple linear extrapolation:
\[
\gamma^L_{J+1} = 2 \: \gamma^L_{J}-\gamma^L_{J-1}, \; \; \; \gamma^R_{J} = 2 \: \gamma^L_{J+1}-\gamma^L_{J+2}.
\]
Using then central differences to discretize (\ref{eqL}) and (\ref{eqR}), linear and quadratic interpolation to discretize respectively the two conditions (\ref{jumpsU_Flux}), we obtain the following second order $(N+4) \times (N+4)$ linear system:
\begin{align}
  u^L_0&=g_0 \label{LS1} \\
  \displaystyle \frac{1}{h^2} \left( \gamma^L_{j- \frac{1}{2}} \left( u^L_j-u^L_{j-1} \right) + \gamma^L_{j+ \frac{1}{2}} \left( u^L_j-u^L_{j+1} \right) \right) &= f_j \; \; \; \; \; \; j=1, \ldots J \label{LS2} \\
  \left( (1-\theta) u^R_J+ \theta u^R_{J+1} \right) - \left((1-\theta) u^L_J+ \theta u^L_{J+1} \right) &= g_D \label{LS3} \\ 
  \gamma^R_\alpha \: \L_{J}' [u^R] (\alpha) - \gamma^L_\alpha \: \L_{J-1}'[u^L](\alpha) &= g_N \label{LS4} \\ 
  \displaystyle \frac{1}{h^2} \left( \gamma^R_{j- \frac{1}{2}} \left( u^R_j-u^R_{j-1} \right) + \gamma^R_{j+ \frac{1}{2}} \left( u^R_j-u^R_{j+1} \right) \right) &= f_j \; \; \; \; \; \; j=J+1, \ldots N \label{LS5} \\
  u^R_N+1 &= g_1, \label{LS6}
\end{align}
with $\gamma^L_\alpha$ and $\gamma^R_\alpha$ obtained by linear interpolation:
\[
\gamma^L_\alpha = (1-\theta) \gamma^L_J+ \theta \gamma^L_{J+1}, \; \; \; \gamma^R_\alpha = (1-\theta) \gamma^R_J+ \theta \gamma^R_{J+1}
\]
and $\theta = (\alpha-x_{J})/h \in [0,1]$.
\\If we apply a simple iterative method such as Gauss-Seidel or Jacobi to this linear system, in general it will not converge, unless we solve the $2 \times 2$ sub-system of transmission conditions, eliminating them from the whole system. This elimination is easy to perform in one dimension, but becomes quite involved in higher dimension. Therefore, we prefer to work with the whole linear system without eliminate transmission conditions from it, in order to extend the method to higher dimension in a forthcoming paper~\cite{CocoRusso:discontinuous2d}. Then we have to find a different approach to solve iteratively the previous linear system. This can be done by relaxing the transmission conditions.

\subsection{Iterative method}\label{ItMethod}
In order to  find a convergent iterative method to solve the linear system (\ref{LS1})-(\ref{LS6}), following the approach introduced in~\cite{CocoRusso:Elliptic} we solve the \textit{associate time-dependent} problem in the unknowns $u^L(x,t)$ and $u^R(x,t)$ for $(x,t) \in [0,1] \times (0,+\infty)$:
\begin{eqnarray}
u^L(0,t)&=&g_0 \label{tDirBC_L}\\
\pad{u^L}{t} &=& \mu \left( \pad{}{x} \left( \gamma^L \pad{u^L}{x} \right) + f \right), \; \; \; x \in [0,\alpha[ \label{tEqL} \\
\left. \pad{u^L}{t} \right|_{x=\alpha} &=& \mu_N \left( \left[ \gamma \pad{u}{x} \right] -g_N\right) \label{tjumpFlux} \\
\left. \pad{u^R}{t} \right|_{x=\alpha}  &=& \mu_D \left( g_D - \left[ u \right] \right) \label{tjumpU} \\
\pad{u^R}{t} &=& \mu \pad{}{x} \left( \gamma^R \pad{u^R}{x} \right) + f, \; \; \; x \in ]\alpha,1] \label{tEqR} \\
u^R(1,t)&=&g_1. \label{tDirBC_R}
\end{eqnarray}
where $\mu$ is a positive function, and $\mu_D$ and $\mu_N$ are two positive constants, that will be set in Sec. \ref{muDmuN} to satisfy some stability condition.
\\The choice of the sign of the two constants $\mu_D$ and $\mu_N$ is crucial and requires some explanation.
Roughly speaking, when replacing a vector equation $F(w) = 0$ for $F:\R^m \rightarrow \R^m$ by $\stackrel{\cdot}{\omega}=F(\omega)$ $\overset{\cdot}{\omega}$, we have to be sure that the solution is asymptotically stable, i.e. that $\lambda(\nabla_\omega F) <0$.
Eq. (\ref{tjumpU}) will be used to compute $u_J^R$, therefore the derivative of the right hand side of Eq. (\ref{tjumpU}) with respect to $u_J^R$ has to be negative, to ensure convergence to equilibrium. Eq. (\ref{tjumpFlux}) is used to determine $u_{J+1}^L$ by a transport equation on $u^L(x,t)$. Since $x_{J+1}>\alpha$ the propagation speed $\mu_N \, \gamma^L$ associated to $u^L(x,t)$, has to be positive.
\\We are obviously interested in the steady-state solution and the time $t$ represents an iterative parameter.
We observe that transmission conditions (\ref{tjumpFlux}) and (\ref{tjumpU}) can be replaced by
\begin{eqnarray*}
\left. \pad{u^R}{t} \right|_{x=\alpha} &=& \mu_N \left( g_N-\left[ \gamma \pad{u}{x} \right] \right) \\
\left. \pad{u^L}{t} \right|_{x=\alpha}  &=& \mu_D \left( \left[ u \right] - g_D \right)
\end{eqnarray*}
because both choices lead to the same steady state conditions.
\\To obtain a second order accurate solution in space we are allowed to discretize first order accurate the time derivative. Using forward Euler in time and central differences in space for (\ref{tEqL}) and (\ref{tEqR}), we obtain (superscripts $L$ and $R$ are omitted):
\begin{equation}\label{dEq1}
u_j^{(m+1)} = u_j^{(m)} + \mu_j \, \Delta t \left( f_j - \displaystyle \frac{\gamma_{j- \frac{1}{2}} \left( u_j^{(m)}-u_{j-1}^{(m)} \right) + \gamma_{j+\frac{1}{2}} \left( u_j^{(m)}-u_{j+1}^{(m)} \right)}{h^2} \right),
\end{equation}
where $j=1,\ldots,J$ for $u^L$ and $j=J+1,\ldots,N$ for $u^R$.
Choosing the maximum time step allowed by the CFL condition for diffusion equation, i.e., $\mu_j \, \Delta t = h^2/(\gamma_{j+1/2}+\gamma_{j-1/2})$, Eq. (\ref{dEq1}) becomes:
\begin{equation}\label{it1}
u_j^{(m+1)} = \displaystyle \frac{1}{\gamma_{j-\frac{1}{2}}+\gamma_{j+ \frac{1}{2}}} \left( f_j \: h^2 +
\gamma_{j- \frac{1}{2}} \: u_{j-1}^{(m)} + \gamma_{j+ \frac{1}{2}} \: u_{j+1}^{(m)} \right),
\end{equation}
where $j=1,\ldots,J$ for $u^L$ and $j=J+1,\ldots,N$ for $u^R$. Observe that such equation is the one obtained by applying Jacobi iteration to Eqs. (\ref{LS2}) and (\ref{LS5}).
\\Let us discretize Eq. (\ref{tjumpFlux}). The time derivative is discretized by forward Euler at the ghost point $x_{J+1}$, which is the quantity we want to compute. The jump is discretized as in (\ref{LS4}), so it is second order accurate. We obtain the iteration:
\begin{equation}\label{it2}
u_{J+1}^{L,(m+1)} = u_{J+1}^{L,(m)} + \mu_N \Delta t \left( \gamma^R_\alpha \: \L_{J}' [u^{R,(m)}] (\alpha) - \gamma^L_\alpha \: \L_{J-1}'[u^{L,(m)}](\alpha) -g_N \right).
\end{equation}
Likewise, in Eq. (\ref{tjumpU}) we discretize the time derivative in $x_J$, obtaining:
\begin{equation}\label{it3}
\begin{split}
u_{J}^{R,(m+1)} &= u_J^{R,(m)} \\
&+ \mu_D \, \Delta t \left( (1-\theta) u_J^{L,(m)} + \theta u_{J+1})^{L,(m)} - (1-\theta) u_J^{R,(m)} + \theta u_{J+1}^{R,(m)} +g_D \right).
\end{split}
\end{equation}
Iterations (\ref{it1}), (\ref{it2}) and (\ref{it3}) constitute the iterative scheme to solve problem (\ref{modprob}) to second order accuracy.

\subsection{Choosing constants $\mu_D$ and $\mu_N$ for transmission conditions}\label{muDmuN}
In (\ref{it2}) and (\ref{it3}) two arbitrary constants $\mu_D$ and $\mu_N$ appear. Following the same argument as in~\cite{CocoRusso:Elliptic}, such constants will be chosen in order to satisfy some stability condition for the equation where they appear. This procedure is not rigorous because it does not take into account the coupling between the equations, and does not consist in a convergence proof. However, in all numerical tests we performed, the conditions we find seem to guarantee convergence.
\\Constant $\mu_D$ is introduced in Eq. (\ref{tjumpU}), which is just a relaxation of the jump condition. Then we require:
\begin{equation}\label{condmuD}
\mu_D \: \Delta t < 1.
\end{equation}
This condition will ensure positivity, and is a factor $2$ more stringent than just stability restriction.
For practical purpose, we set $\mu_D \: \Delta t = 0.9$.
In order to obtain a condition on $\mu_N$, we rewrite Eq. (\ref{tjumpFlux}) as follows (we have supposed for simplicity homogeneous jump $g_N=0$):
\begin{equation}\label{conv}
\pad{u^L}{t} + \mu_N \: \gamma^L \: \pad{u^L}{x} = \mu_N \: \gamma^R \: \pad{u^R}{x}, \; \; \; t \in (0,\infty).
\end{equation}
This is a simple convection equation with speed $\mu_N \: \gamma^L$. Then a simple CFL condition for convection equation might be 
\[
\mu_N \Delta t \leq \displaystyle \frac{h}{\gamma^L}.
\]
Numerical experiments show that this condition is not enough, especially in the case $\gamma^R / \gamma^L \gg 1$. An explanation of this behavior may be that the right-hand side of (\ref{conv}) is not stationary when the convection evolves in time, but it depends on time itself by $u^R$. An acceptable condition is 
\begin{equation}\label{condmuN}
\mu_N \, \Delta t \leq \displaystyle \frac{h}{\max \left\{ \gamma^L, \gamma^R \right\} }.
\end{equation}
For practical purpose we choose $\mu_N \, \Delta t =0.9 \, h / \max \left\{ \gamma^L, \gamma^R \right\}$.
Numerical tests show that conditions (\ref{condmuD}) and (\ref{condmuN}) are sufficient for guarantee convergence, but not necessary.
A more detailed analysis is in progress.
\\Notice that $\mu \, \Delta t = O(h^2)$, $\mu_N \, \Delta t=O(h)$, $\mu_D \, \Delta t=O(1)$. Furthermore, only the product of the constants times $\Delta t$ enters into the conditions, therefore we may imagine that $\Delta t=1$.

\section{Multigrid approach}
The convergence of the iterative method proposed in Sec. \ref{ItMethod} is usually very slow. To accelerate the convergence we use a multigrid strategy. To make the iteration scheme (\ref{it1})-(\ref{it3}) a building block for an efficient multi-grid solver, we must be sure that such iteration (\textit{relaxation scheme}) has the \textit{smoothing property}, i.e. that after few steps, the error becomes smooth (not necessarily small). Roughly speaking, the high-frequency components of the error reduce quickly. We do not explain all multigrid features, but just what is different from classical multigrid approach, remanding to the literature for more details (e.g., see~\cite{Trottemberg:MG, Hackbusch:MG, Briggs:MG}). The iteration scheme (\ref{it1})-(\ref{it3}) is a Jacobi-like scheme, as mentioned in Sec. \ref{ItMethod}. Jacobi scheme is not a good smoother, since high-frequency components of the error reduce slowly. A good smoother is instead the Gauss-Seidel scheme. Then, we use a Gauss-Seidel version of (\ref{it1})-(\ref{it3}) as relaxation scheme, i.e.
\begin{align}
u_j^{L,(m+1)} &= \displaystyle \frac{1}{\gamma_{j-\frac{1}{2}}+\gamma_{j+ \frac{1}{2}}} \left( f_j \: h^2 +
\gamma_{j- \frac{1}{2}} \: u_{j-1}^{L,(m+1)} + \gamma_{j+\frac{1}{2}} \: u_{j+1}^{L,(m)} \right), \; \; \; j=1,\ldots,J \label{it1GS} \\
u_{J+1}^{L,(m+1)} &= u_{J+1}^{L,(m)} + \mu_N \Delta t \left( \gamma^R_\alpha \: \L_{J}' [u^{R,(m)}] (\alpha) - \gamma^L_\alpha \: \L_{J-1}'[\tilde{u}^{L}](\alpha) -g_N \right) \label{it2GS} \\
\begin{split}
u_{J}^{R,(m+1)} &=u_J^{R,(m)} \\
&+ \mu_D \Delta t \left( (1-\theta) u_J^{L,(m+1)} + \theta u_{J+1}^{L,(m+1)} - (1-\theta) u_J^{R,(m)} - \theta u_{J+1}^{R,(m)} +g_D \right) 
\end{split}\label{it3GS} \\
u_j^{R,(m+1)} &= \displaystyle \frac{1}{\gamma_{j-\frac{1}{2}}+\gamma_{j+\frac{1}{2}}} \left( f_j \: h^2 +
\gamma_{j-\frac{1}{2}} \: u_{j-1}^{R,(m+1)} + \gamma_{j+\frac{1}{2}} \: u_{j+1}^{R,(m)} \right), \; \; \; j=J+1,\ldots,N \label{it4GS}
\end{align}
where in (\ref{it2GS}) we intend $\tilde{u}^L$ such that $\tilde{u}^L_j = u_j^{L,(m+1)}$ for $j < J+1$ and $\tilde{u}^L_{J+1}=u^{L,(m)}_{J+1}$.
The unknowns are updated in the same order reported in (\ref{unknowns}).
\\In order to explain the multigrid approach, we just describe the two-grid correction scheme (TGCS), because all the other schemes, such as $V$-cycle, $W$-cycle, $F$-cycle or Full Multigrid cycle, can be easily derived from it (see~\cite[Sections 2.4 and 2.6]{Trottemberg:MG} for more details). Let us introduce some notation. For a grid of spatial step $h$, we denote: 
\[
J=\left\lfloor \displaystyle \frac{\alpha}{h} \right\rfloor, \; \; \; \theta=\displaystyle \frac{\alpha}{h}-J
\]
\[
S(\Omega_h) = \left\{ \wh = (w^L, w^R) \mbox{ such that }  w^L \colon \left\{ x_0, \ldots, x_{J+1} \right\} \rightarrow \R, \; w^R \colon \left\{ x_J, \ldots, x_{N+1} \right\} \rightarrow \R   \right\}
\]
\[
\stackrel{\circ}{S}(\Omega_h) = \left\{ \wh = (w^L, w^R) \mbox{ such that }  w^L \colon \left\{ x_1, \ldots, x_{J} \right\} \rightarrow \R, \; w^R \colon \left\{ x_{J+1}, \ldots, x_{N} \right\} \rightarrow \R   \right\}
\]
\[
\uh = ( (u^L_j)_{j=0,\ldots,J+1}, (u^R_j)_{j=J,\ldots,N+1}) \in S(\Omega_h)
\]
\[
\gamma_h = ( (\gamma^L_j)_{j=0,\ldots,J+1}, (\gamma^R_j)_{j=J,\ldots,N+1}) \in S(\Omega_h)
\]
\[
\fh \in \stackrel{\circ}{S}(\Omega_h) \mbox{ such that } \fh(x_j)=f_j
\]
\[
L_h \colon S(\Omega_h) \times S(\Omega_h) \longrightarrow \stackrel{\circ}{S}(\Omega_h) \mbox{ such that }
\]
\begin{eqnarray*}
\left( L_h (\gamma_h, \uh) \right)_j &=& \displaystyle \frac{1}{h^2} \left( \gamma^{L}_{j- \frac{1}{2}} \left( u_j^{L}-u_{j-1}^{L} \right) + \gamma^{L}_{j+\frac{1}{2}} \left( u_j^{L}-u_{j+1}^{L} \right) \right) \mbox{ if } j \leq J \\
\left( L_h (\gamma_h, \uh) \right)_j &=& \displaystyle \frac{1}{h^2} \left( \gamma^{R}_{j- \frac{1}{2}} \left( u_j^{R}-u_{j-1}^{R} \right) + \gamma^{R}_{j+ \frac{1}{2}} \left( u_j^{R}-u_{j+1}^{R} \right) \right) \mbox{ if } j \geq J+1
\end{eqnarray*}
\[
\left[ \: \cdot \: \right]^D_h \colon S(\Omega_h) \longrightarrow \R \mbox{ such that }
\]
\[
\left[ \uh \right]^D_h = \left( (1-\theta) u^R_J+ \theta u^R_{J+1} \right) - \left((1-\theta) u^L_J+ \theta u^L_{J+1} \right)
\]
\[
\left[ \: \cdot \: , \: \cdot \: \right]^N_h \colon S(\Omega_h) \times S(\Omega_h) \longrightarrow \R \mbox{ such that }
\]
\[
\left[ \gamma_h, \uh \right]^N_h = \gamma^R_\alpha \: \L_{J}' [u^R] (\alpha) - \gamma^L_\alpha \: \L_{J-1}'[u^L](\alpha)
\]
The linear system (\ref{LS1})-(\ref{LS6}) can be resumed as follows:
\begin{eqnarray}
L_h (\gamma_h, \uh) &=& \fh \\
\left[ \uh \right]^D_h &=& g_D \\
\left[ \gamma_h , \uh \right]^N_h &=& g_N \\
u^L_0&=&g_0 \\
u^R_N&=&g_1.
\end{eqnarray}
For simplicity we assume that $N+1=1/h$ is a power of $2$.
The TGCS consists into the following algorithm:
\begin{enumerate}\label{alg1d}
\item Set initial guess $\uh = 0$. \\
\item Relax $\nu_1$ times on the finest grid:
for $k$ from $1$ to $\nu_1$ do (\ref{it1GS}), (\ref{it2GS}), (\ref{it3GS}).
\item Compute the defects $\rh \in \stackrel{\circ}{S}(\Omega_h)$, $\tilde{g}_D, \tilde{g}_N \in \R$: 
\begin{eqnarray*}
\rh &=& \fh + L_h(\gamma_h,\uh) \\
\tilde{g}_D &=& g_D - \left[ \uh \right]^D_h \\
\tilde{g}_N &=& g_N - \left[ \gamma_h , \uh \right]^N_h
\end{eqnarray*}
\item Transfer the defect $\rh$ to a coarser grid with spatial step $2h$ by a suitable \textit{restriction operator}
\[
\rhh = I_{2h}^{h} \left( \rh \right).
\]
\item Solve exactly the residual problem on the coarser grid in the unknow $\ehh \in S(\Omega_{2h})$
\begin{eqnarray*}
L_h (\gamma_{2h}, \ehh) &=& \rhh \\
\left[ \ehh \right]^D_h &=& \tilde{g}_D \\
\left[ \gamma_{2h} , \uhh \right]^N_h &=& \tilde{g}_N \\
e^L_0&=&0 \\
e^R_{(N+1)/2}&=&0
\end{eqnarray*}
\item Transfer the error to the finest grid by a suitable \textit{interpolation operator}
\[
\eh = I_h^{2h} \left( \ehh \right).
\] 
\item Correct the fine-grid approximation
\[
\uh = \uh + \eh.
\]
\item Relax $\nu_2$ times on the finest grid:
for $k$ from $1$ to $\nu_2$ do (\ref{it1GS}), (\ref{it2GS}), (\ref{it3GS}).
\end{enumerate}
To complete the description of TGCS, we have just to explain the steps concerning grid migration (steps 4 and 6).

\subsection{Transfer grid operators}
In this section, we describe the transfer grid operators for vertex-centered grid. Observe that coefficients $\gamma^L$ and $\gamma^R$ can be transferred in an exact manner by a simple injection operator.

\subsubsection{Restriction operator}
\begin{figure}[!ht]
 \begin{minipage}[b]{0.99\textwidth}
   	\centering
		\includegraphics[width=0.60\textwidth]{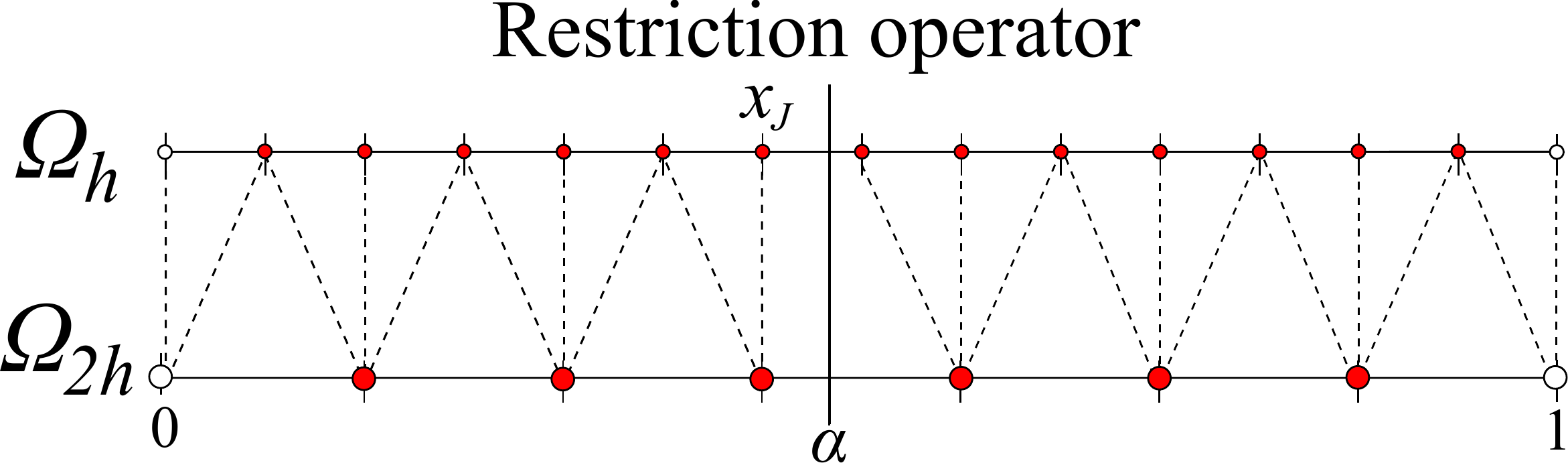}
 \end{minipage}
 \vskip 0.5 cm
  \begin{minipage}[b]{0.99\textwidth}
   	\centering
		\includegraphics[width=0.60\textwidth]{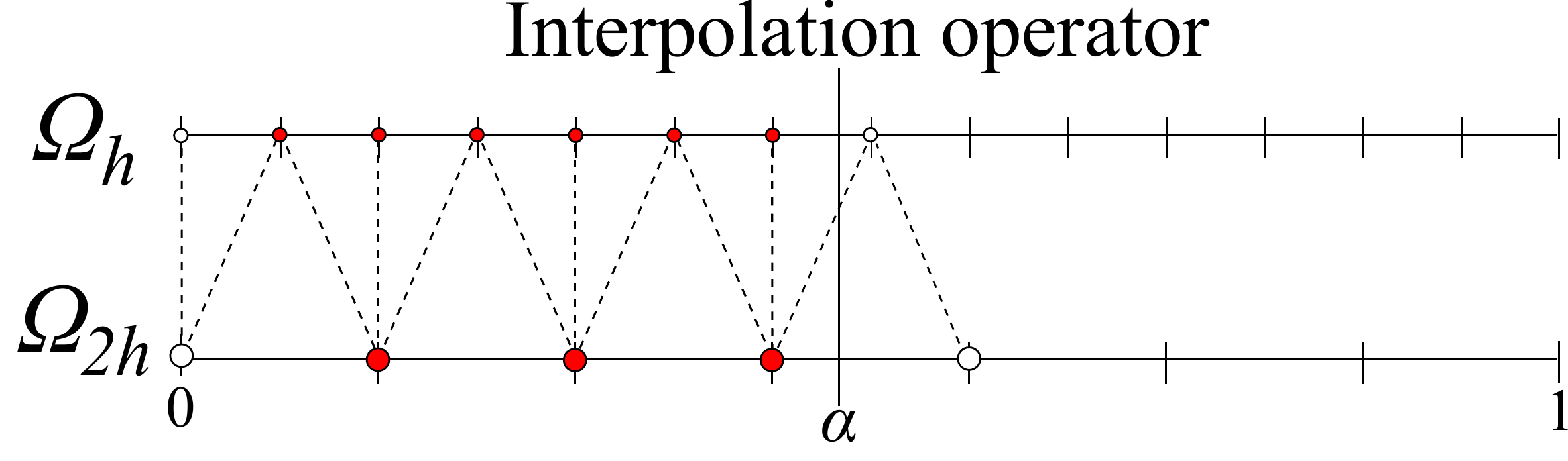}
 \end{minipage}
 \vskip 0.5 cm
 \begin{minipage}[b]{0.99\textwidth}
  	\centering
		\includegraphics[width=0.60\textwidth]{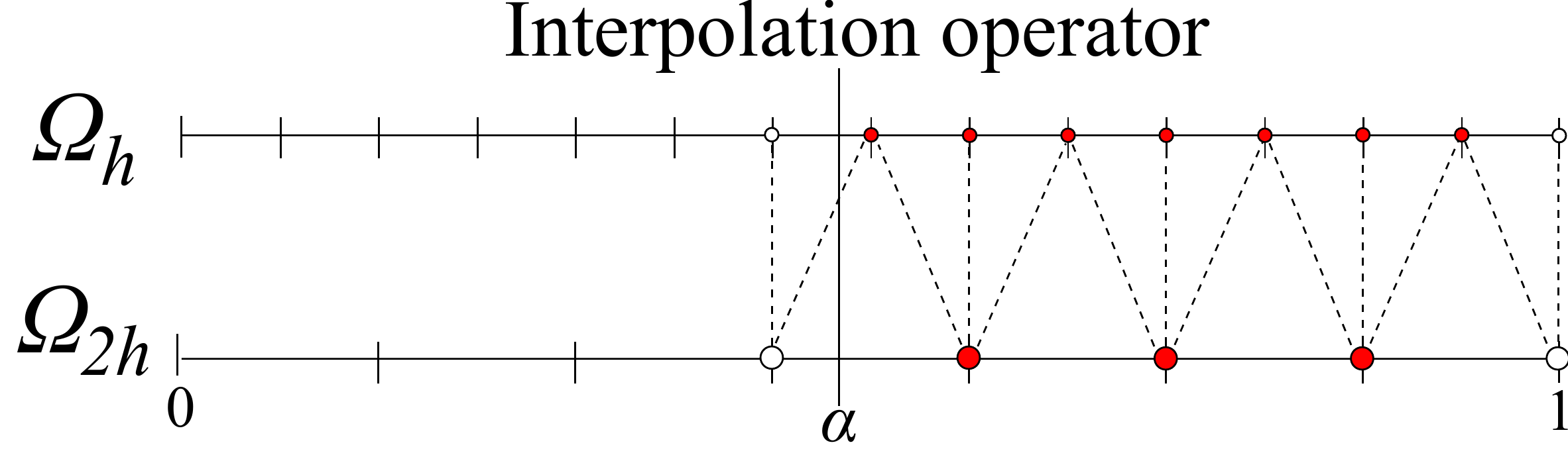}
	\caption{\footnotesize{ Fine and coarse grid for transfer operators. The dashed lines represent the action of the restriction (top) and the interpolation (middle and bottom) operators. } }
	\label{fig:VC}
 \end{minipage}
\end{figure}

Since such operator will act on the defect $\rh=(\rh^L,\rh^R) \in \stackrel{\circ}{S}(\Omega_h)$ (step 4), we perform the restriction from a fine grid to a coarser grid separately for $\rh^L$ and $\rh^R$. This is justified by the fact that the defect $\rh^L$ of the left domain may be very different (after few relaxations) from the defect $\rh^R$ of the right domain, especially in the case of high jumping coefficient, i.e., $\max \left\{ \gamma_\alpha^L / \gamma_\alpha^R, \gamma_\alpha^R / \gamma_\alpha^L \right\} >> 1$. In addition, these defects are very different also from the defects of jumping conditions $\tilde{g}_D$ and $\tilde{g}_N$, because the operators scale with different power of $h$.
\\Let us describe the restriction of $\rh^L$ by the operator $\left(I^h_{2h} \right)^L$ (see Fig. \ref{fig:VC}).
Let $x_J$ be the closest grid point to $\alpha$ from the left in the fine grid (see Fig. \ref{fig:Omega1d}). Let $x$ be a grid point of the coarse grid. If $x < x_J$ we will use the standard full-weighting restriction operator (FW):
\begin{equation}\label{restL}
\left(I^h_{2h} \right)^L \rh^L (x) = \displaystyle \frac{1}{4} \left( \rh(x-h)^L + 2 \, \rh(x)^L + \rh(x+h)^L \right),
\end{equation}
while if $x=x_J$ we reduce to an \textit{upwind} linear convex combination from the left direction:
\begin{equation}\label{restredL}
\left(I^h_{2h} \right)^L \rh^L (x) = \omega_1 \, \rh^L(x) + (1-\omega_1) \rh^L(x-h),
\end{equation}
since in $x+h$ only $\rh^R$ is defined and not $\rh^L$. In our tests we found that $\omega_1=1/2$ gives better results than $\omega_1=3/4$.
\\The operator $\left(I^h_{2h} \right)^R$ works in a similar manner: let $x_{J+1}$ the closest grid point to $\alpha$ from the right in the fine grid. If $x>x_{J+1}$ we will use the standard full-weighting restriction operator (FW):
\begin{equation}\label{restR}
\left(I^h_{2h} \right)^R \rh^R (x) = \displaystyle \frac{1}{4} \left( \rh(x-h)^R + 2 \, \rh(x)^R + \rh(x+h)^R \right),
\end{equation}
while if $x=x_J$ we reduce to an \textit{Upwind} mean value from the left direction:
\begin{equation}\label{restredR}
\left(I^h_{2h} \right)^R \rh^R (x) =  \frac{1}{2} \left( \rh^R(x) + \rh^R(x+h) \right).
\end{equation}
The whole restriction reads
\[
I^h_{2h} \rh = \left( \left(I^h_{2h} \right)^L \rh^L, \left(I^h_{2h} \right)^R \rh^R \right).
\]
In the upper part of Fig. \ref{fig:VC} is represented the case in which we have to use (\ref{restredL}) and (\ref{restR}). The only other possible case is that we have to use (\ref{restL}) and (\ref{restredR}).

\subsubsection{Interpolation operator}
Since such operator will act on the correction $\ehh=(\ehh^L,\ehh^R) \in S(\Omega_{2h})$ (step 4), we perform the interpolation from a coarse grid to a finer grid separately for $\ehh^L$ and $\ehh^R$ (see middle and lower part of Fig. \ref{fig:VC}), but always using the standard linear interpolation:

\begin{equation}\label{interpL}
\left\{
\begin{array}{rcll}
\left( I^{2h}_{h}\right)^L \ehh^L (x_j) &=& \ehh^L(x_j) & \mbox{ if $j$ is even} \\
\left( I^{2h}_{h} \right)^L \ehh^L (x_j) &=&  \frac{1}{2} \left( \ehh^L(x_{j-1}) + \ehh^L(x_{j+1}) \right) & \mbox{ if $j$ is odd}. 
\end{array}
\right.
\end{equation}

\begin{equation}\label{interpR}
\left\{
\begin{array}{rcll}
\left( I^{2h}_{h}\right)^R \ehh^R (x_j) &=& \ehh^R(x_j) & \mbox{ if $j$ is even} \\
\left( I^{2h}_{h} \right)^R \ehh^R (x_j) &=&  \frac{1}{2} \left( \ehh^R(x_{j-1}) + \ehh^R(x_{j+1}) \right) & \mbox{ if $j$ is odd}. 
\end{array}
\right.
\end{equation}

The whole interpolation reads
\[
I^{2h}_{h} \ehh = \left( \left(I_h^{2h} \right)^L \ehh^L, \left(I_h^{2h} \right)^R \ehh^R \right).
\]

\textbf{Remark. 1 (Coarser operator)} We observe that the discrete operator $L_{2h}$ on the coarser grid (step 5) is just the operator obtained discretizing directly the continuous operator in the grid with spatial step $2h$, and not the operator obtained by Galerkin condition
\begin{equation}\label{Galerkin}
L_{2h} = I^{h}_{2h} \: L_{h} \: I^{2h}_{h}.
\end{equation}
The last approach, typical of algebraic multigrid, makes the algebraic problem more expensive from a computational point of view and does not take advantage of the fact that the discrete problem comes from a continuous problem.

\textbf{Remark. 2 ($V$-cycle)} The $V$-cycle algorithm is easily obtained from the TGCS recursively, namely applying the same algorithm to solve the residual equation in step 5. To terminate the recursion, an exact solver is used to solve the residual problem when the grid achieves a fixed level of coarsening. We denote by $V(\nu_1,\nu_2)$-cycle the $V$-cycle performed with $\nu_1$ pre-relaxations and $\nu_2$ post-relaxations.  

\textbf{Remark. 3 ($W$-cycle)} The $W$-cycle is similar to the $V$-cycle, with the only difference that the residual problem is solved recursively two times instead of one (in general schemes, $\delta$ times, but $\delta > 2$ is considered useless for practical purposes).

\section{Numerical tests}\label{NumTests}
In this section we confirm numerically the second order accuracy of the discretization of Sec. \ref{discr} and compute the convergence factor $\rho$ of the multigrid approach for several examples, to confirm the independence of $\rho$ from the spatial step $h$ and the magnitude of the jumping coefficient.
\\Second order accuracy is gained also for first derivative of the solution, as it is shown by the comparison between exact first derivative and the numerical derivative obtained by central difference of the numerical solution.
\\In all numerical tests, we choose an arbitrary interface $\alpha \in ]0,1[$ and an analytical expression of the exact solution $u=(u^L,u^R)$ and of diffusion coefficient $\gamma=(\gamma^L,\gamma^R)$. Then we reconstruct the data $f$, $g_D$ and $g_N$, perform the multigrid technique, and compare the numerical solution with the exact solution to compute the order of accuracy by the slope of the best-fit line. 
In all our tests we use the following stopping criterion for the $V-$cycle
\[
\displaystyle \frac{ \left\| \uh^{(m+1)}-\uh^{(m)} \right\|_\infty }{ \left\| \uh^{(m+1)} \right\|_\infty } \leq TOL.
\]
This will ensure that the actual relative error satisfies
\[
\displaystyle \frac{ \left\| \eh^{(m+1)} \right\|_\infty }{ \left\| \eh \right\|_\infty } \leq \rho \displaystyle \frac{TOL}{1-\rho}.
\]
The tolerance we used is $TOL=10^{-6}$, which ensures that the error in the solution of the algebraic system is always lower than truncation error.
For each example we show a table in which we list the errors, and the value in the third [fifth] column and $i$-th row of the table indicates the accuracy order, computed as $\log_2 \left( e_{i-1} /e_{i}  \right)$, where $e_i$ is the $L^\infty$-error of the numerical solution [derivative] indicated in the second [fourth] column and $i$-th row.
\\To compute the asymptotic convergence factor, we use the following estimate:
\[
\rho=\rho^{(m)} = \displaystyle \frac{\left\| \rh^{(m)} \right\|_{\infty}}{\left\| \rh^{(m-1)} \right\|_{\infty}},
\]
which is reliable for $m$ large.
In order to avoid difficulties related to numerical instability due to machine precision, we will always use the homogeneous model problem as a test when we want to compute the asymptotic convergence factor, namely Eq. (\ref{modprob}) with $f=g_0=g_1=0$ and homogeneous jump conditions, and perform the multigrid algorithm starting from an initial guess different from zero. Since in this case we are just interested in the convergence factor and not in the numerical solution itself (which approaches zero), a reasonable stop criterion will be
\[
\displaystyle \frac{ \left|\rho^{(m)}-\rho^{(m-1)} \right| }{ \rho^{(m)} } < 10^{-2}.
\]
Several tests are performed for each example, based on the different size of the finest and coarsest grids. The finest grid is obtained dividing the domain $[0,1]$ into $N+1$ intervals, while the coarsest grid is obtained dividing the domain into $N_c+1$ intervals.

\subsection{Example 1}\label{ex1}
\begin{figure}[!hbt]
\begin{minipage}{0.49\textwidth}
   	\centering
   	\includegraphics[width=1.00\textwidth]{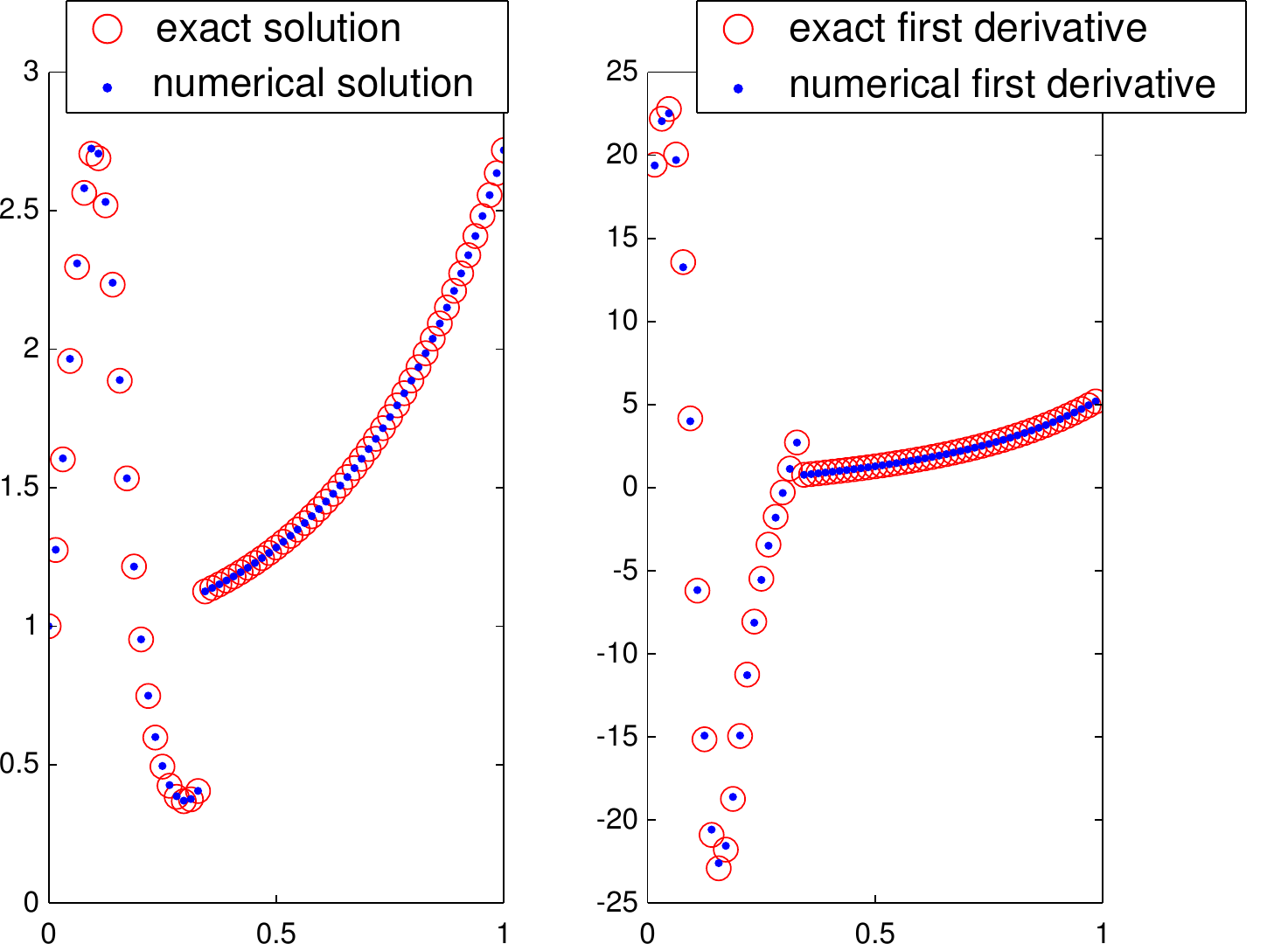}
 \end{minipage}
 \ \hspace{2mm} 
 \begin{minipage}{0.49\textwidth}
  	\centering
  	\captionsetup{width=0.80\textwidth}
	\includegraphics[width=1.00\textwidth]{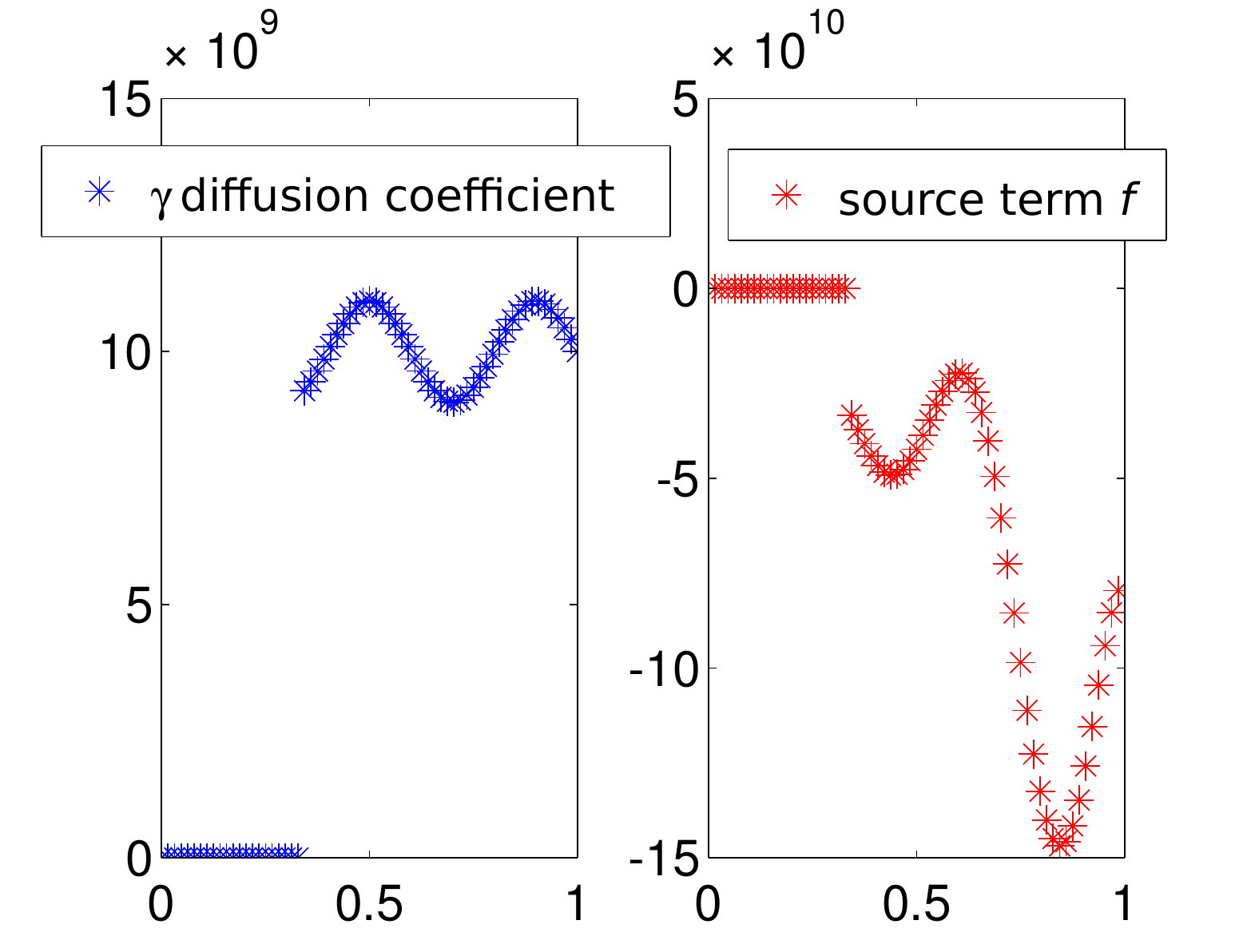}
 \end{minipage}
 \captionsetup{width=0.80\textwidth}
 \caption{\footnotesize{We refer to Ex. \ref{ex1}. The data are computed for $N=64$.}}
 	\label{fig:dataEx1}
\end{figure}

\begin{figure}[!hbt]
 \captionsetup{width=0.80\textwidth}
 \begin{minipage}{0.45\textwidth}
 \centering
   	\includegraphics[width=1.00\textwidth]{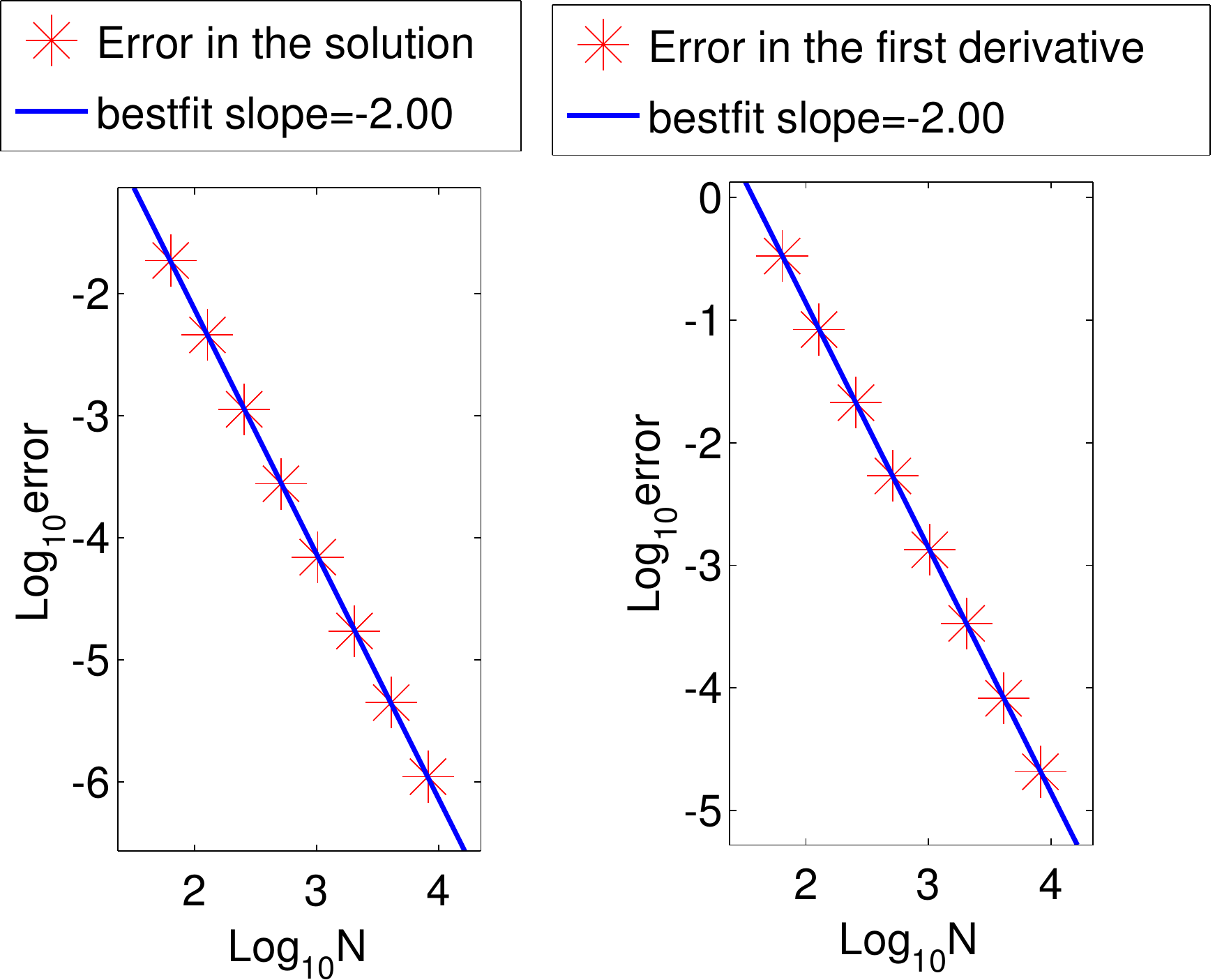}
 \end{minipage}
   \hspace{2mm} 
 \begin{minipage}{0.45\textwidth}
 \centering
\begin{tabular}{||c || c | c || c | c ||}  
\hline\hline                        
$N+1$ & $\left\| \vec{u}-\uh \right\|_\infty$ & order & $\left\| \vec{u}'-\uh' \right\|_\infty$ & order \\ [0.5ex] 
\hline                    
64 & 1.87 $\cdot 10^{-2}$ & - & 3.33 $\cdot 10^{-1}$  & -   \\  
128 & 4.59 $\cdot 10^{-3}$ & 2.03 & 8.38 $\cdot 10^{-2}$  & 1.99   \\  
256 & 1.13 $\cdot 10^{-3}$ & 2.02 & 2.12 $\cdot 10^{-2}$  & 1.98   \\  
512 & 2.77 $\cdot 10^{-4}$ & 2.03 & 5.35 $\cdot 10^{-3}$  & 1.99   \\  
1024 & 6.95 $\cdot 10^{-5}$ & 2.00 & 1.34 $\cdot 10^{-3}$  & 2.00   \\  
2048 & 1.73 $\cdot 10^{-5}$ & 2.01 & 3.36 $\cdot 10^{-4}$  & 1.99   \\  
4096 & 4.48 $\cdot 10^{-6}$ & 1.95 & 8.21 $\cdot 10^{-5}$  & 2.03   \\  
8192 & 1.11 $\cdot 10^{-6}$ & 2.01 & 2.07 $\cdot 10^{-5}$  & 1.99   \\  
\hline     
 \end{tabular} 
 \end{minipage}
	\caption{\footnotesize{We refer to Ex \ref{ex1}.
	Left: Representation of the $L^\infty$-error of the numerical solution and its derivative. The slope of the best-fit lines is respectively $s=-2.00$ and $s=-2.00$.
	 Right: List of errors and order of accuracy computed by subsequent errors. 
	}}
	\label{fig:bestfitEx1}
 	\end{figure}

We choose (see Fig. \ref{fig:dataEx1})
\[
\alpha=0.343, \; \; \;
\left\{
\begin{array}{ccc}
u^L&=&e^{\sin(5 \pi x)} \\
u^R&=&e^{x^2}
\end{array},
\right. \; \; \;
\left\{
\begin{array}{ccc}
\gamma^L&=&3+\cos(5 \pi x) \\
\gamma^R&=&10^9\left( 10+\sin(5 \pi x) \right)
\end{array}.
\right.
\]
Fig. \ref{fig:bestfitEx1} shows the numerical results and the second order slope of the best-fit line for the $L^\infty$-error of the numerical solution and its derivative. 
Table \ref{table:rhoEx1} shows the convergence factor for different values of $N$ and $N_c$.

\begin{table}[!hbt]
\captionsetup{width=0.80\textwidth}
\caption{ \footnotesize{ Measured $V(1,1)$-cycle convergence factor for the numerical test of Ex. \ref{ex1}. We use $N+2$ number of grid points in the finest grid; $N_c+2$ number of grid points in the coarsest grid.  }} 
\centering      
\begin{tabular}{|| c c || c | c | c | c | c | c | c | c ||}  
\hline\hline                        
 & N+1 & 32 & 64 & 128 & 256 & 512 & 1024 & 2048 & 4096  \\ [0.5ex] 
Nc+1 & &  &  &  &  &  &  &  & \\ 
\hline\hline                        
16 & & 0.15 & 0.16 & 0.15 & 0.17 & 0.19 & 0.15 & 0.15 & 0.15\\ 
32 & &  & 0.14 & 0.12 & 0.17 & 0.19 & 0.15 & 0.15 & 0.15\\ 
64 & &  &  & 0.07 & 0.16 & 0.19 & 0.15 & 0.15 & 0.15\\ 
128 & &  &  &  & 0.11 & 0.17 & 0.15 & 0.15 & 0.15\\ 
\hline\hline                        
 \end{tabular} 
      \label{table:rhoEx1}  
 \end{table}

\subsection{Example 2}\label{ex2}
\begin{figure}[!hbt]
\begin{minipage}{0.49\textwidth}
   	\centering
   	\includegraphics[width=1.00\textwidth]{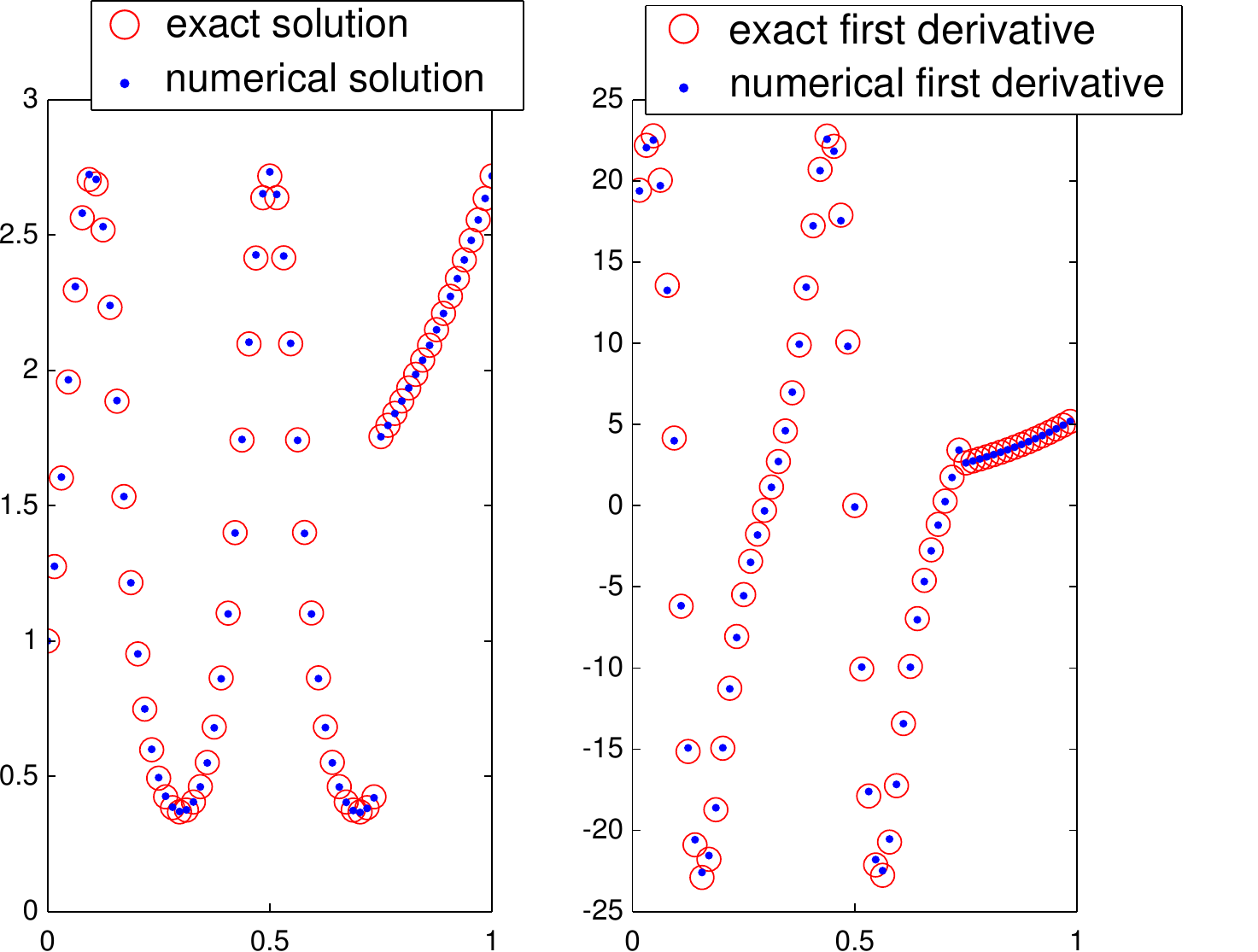}
 \end{minipage}
 \ \hspace{2mm} 
 \begin{minipage}{0.49\textwidth}
  	\centering
  	\captionsetup{width=0.80\textwidth}
	\includegraphics[width=1.00\textwidth]{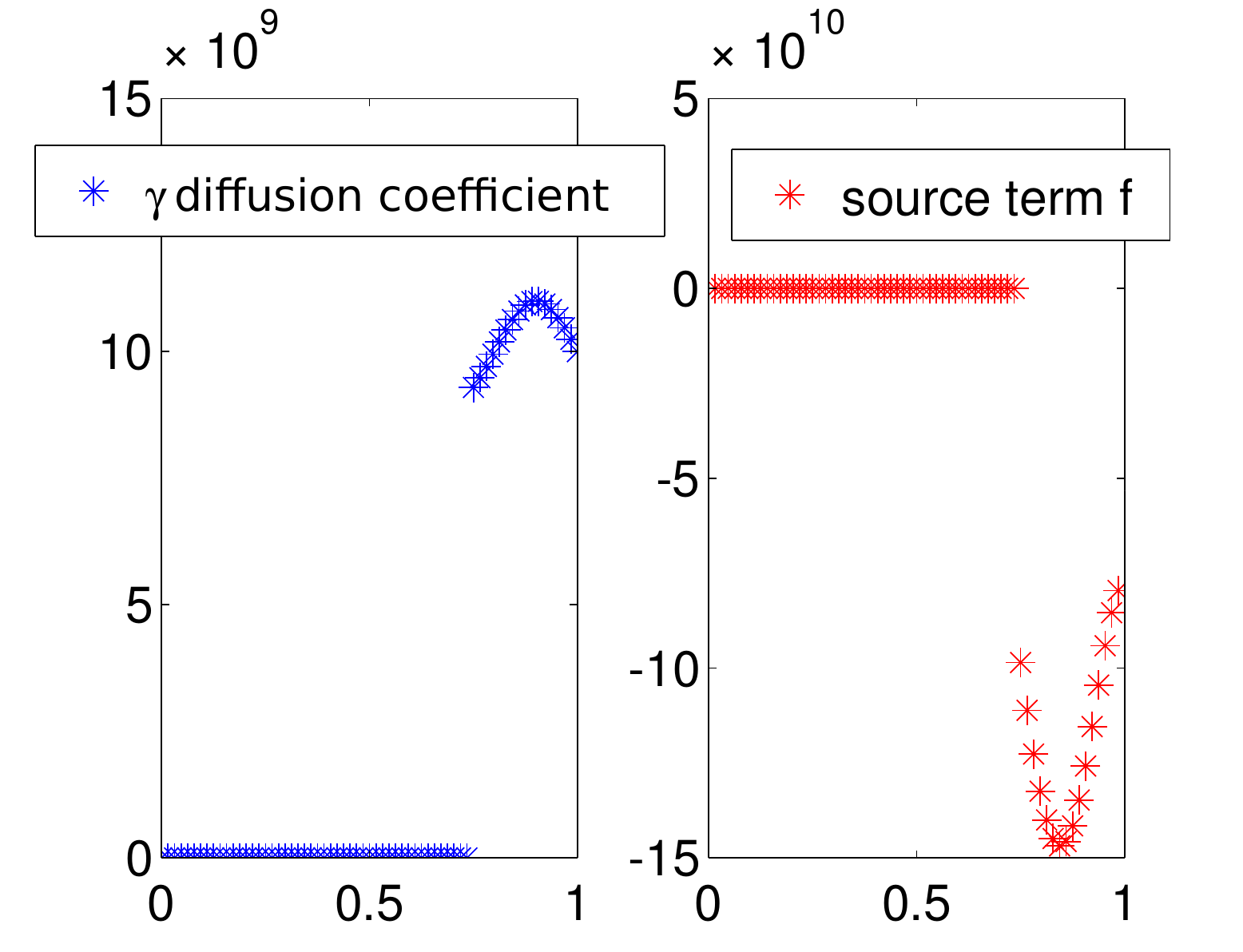}
 \end{minipage}
 \captionsetup{width=0.80\textwidth}
 \caption{\footnotesize{We refer to Ex. \ref{ex2}. The data are computed for $N=64$.}}
 	\label{fig:dataEx2}
\end{figure}

\begin{figure}[!hbt]
 \captionsetup{width=0.80\textwidth}
 \begin{minipage}{0.45\textwidth}
 \centering
   	\includegraphics[width=1.00\textwidth]{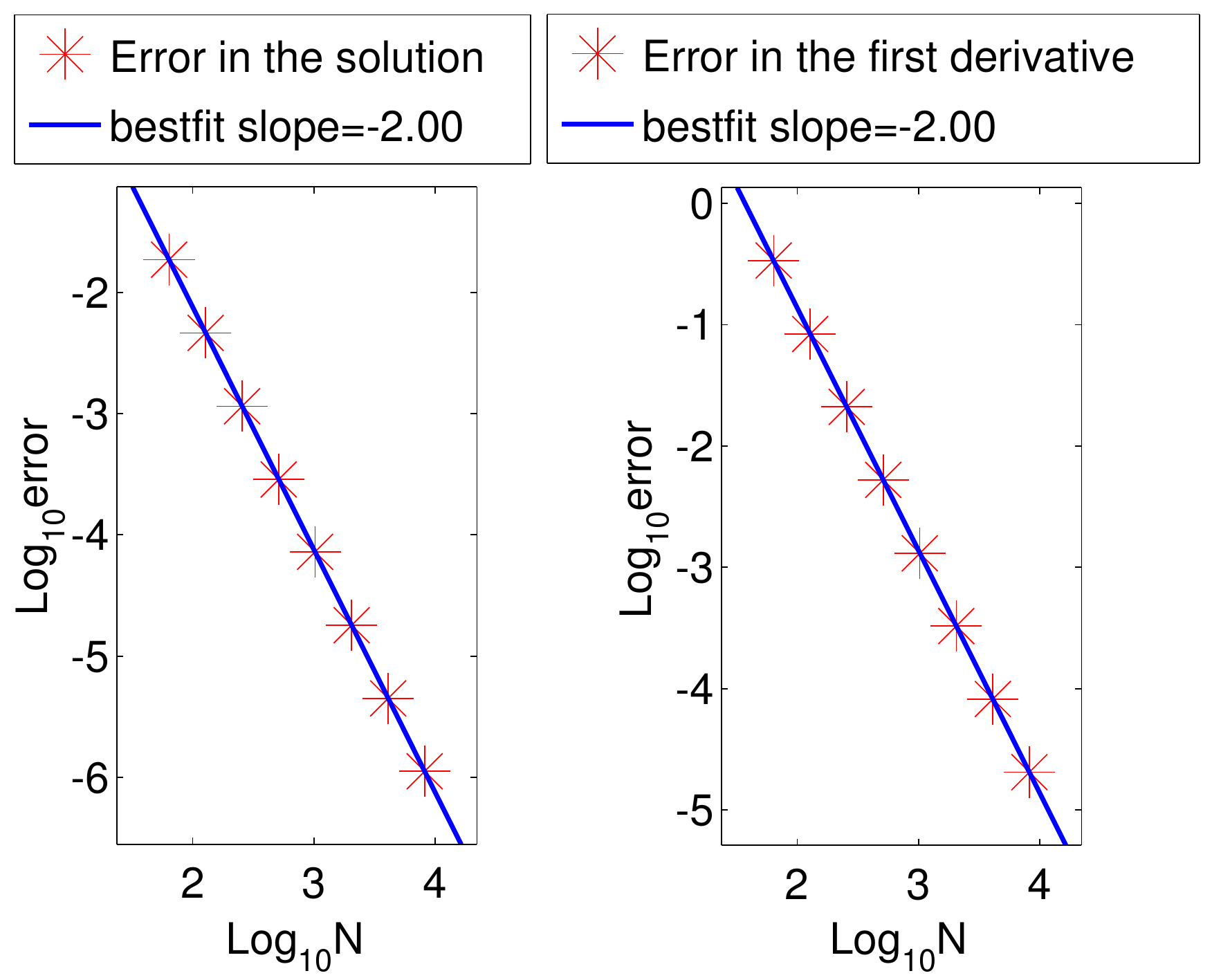}
 \end{minipage}
   \hspace{2mm} 
 \begin{minipage}{0.45\textwidth}
 \centering
\begin{tabular}{||c || c | c || c | c ||}  
\hline\hline                        
$N+1$ & $\left\| \vec{u}-\uh \right\|_\infty$ & order & $\left\| \vec{u}'-\uh' \right\|_\infty$ & order \\ [0.5ex] 
\hline                    
64 & 1.86 $\cdot 10^{-2}$ & - & 3.38 $\cdot 10^{-1}$  & -   \\  
128 & 4.63 $\cdot 10^{-3}$ & 2.01 & 8.38 $\cdot 10^{-2}$  & 2.01   \\  
256 & 1.15 $\cdot 10^{-3}$ & 2.01 & 2.10 $\cdot 10^{-2}$  & 2.00   \\  
512 & 2.86 $\cdot 10^{-4}$ & 2.01 & 5.26 $\cdot 10^{-3}$  & 2.00   \\  
1024 & 7.24 $\cdot 10^{-5}$ & 1.98 & 1.30 $\cdot 10^{-3}$  & 2.01   \\  
2048 & 1.80 $\cdot 10^{-5}$ & 2.01 & 3.28 $\cdot 10^{-4}$  & 1.99   \\  
4096 & 4.48 $\cdot 10^{-6}$ & 2.01 & 8.21 $\cdot 10^{-5}$  & 2.00   \\  
8192 & 1.12 $\cdot 10^{-6}$ & 2.00 & 2.05 $\cdot 10^{-5}$  & 2.00   \\  
\hline     
 \end{tabular} 
 \end{minipage}
	\caption{\footnotesize{We refer to Ex \ref{ex2}.
	Left: Representation of the $L^\infty$-error of the numerical solution and its derivative. The slope of the best-fit lines is respectively $s=-2.00$ and $s=-2.00$.
	 Right: List of errors and order of accuracy computed by subsequent errors. 
	}}
	\label{fig:bestfitEx2}
 	\end{figure}

We choose (see Fig. \ref{fig:dataEx2})
\[
\alpha=0.743, \; \; \;
\left\{
\begin{array}{ccc}
u^L&=&e^{\sin(5 \pi x)} \\
u^R&=&e^{x^2}
\end{array},
\right. \; \; \;
\left\{
\begin{array}{ccc}
\gamma^L&=&3+\cos(5 \pi x) \\
\gamma^R&=&10^9\left( 10+\sin(5 \pi x) \right)
\end{array}.
\right.
\]
The only difference with respect to the previous example is the value of $\alpha$.
\\Fig. \ref{fig:bestfitEx2} shows the numerical results and the second order slope of the best-fit line for the $L^\infty$-error of the numerical solution and its derivative. 
Table \ref{table:rhoEx2} shows the convergence factor for different values of $N$ and $N_c$.

\begin{table}[!hbt]
\captionsetup{width=0.80\textwidth}
\caption{ \footnotesize{ Measured $V(1,1)$-cycle convergence factor for the numerical test of Ex. \ref{ex2}. We use $N+2$ number of grid points in the finest grid; $N_c+2$ number of grid points in the coarsest grid.  }} 
\centering      
\begin{tabular}{|| c c || c | c | c | c | c | c | c | c ||}  
\hline\hline                        
 & N+1 & 32 & 64 & 128 & 256 & 512 & 1024 & 2048 & 4096  \\ [0.5ex] 
Nc+1 & &  &  &  &  &  &  &  & \\ 
\hline\hline                        
16 & & 0.13 & 0.13 & 0.11 & 0.14 & 0.15 & 0.15 & 0.15 & 0.15\\ 
32 & &  & 0.13 & 0.11 & 0.15 & 0.15 & 0.15 & 0.15 & 0.15\\ 
64 & &  &  & 0.11 & 0.13 & 0.15 & 0.15 & 0.15 & 0.15\\ 
128 & &  &  &  & 0.15 & 0.15 & 0.15 & 0.15 & 0.15\\ 
\hline\hline                        
 \end{tabular} 
    \label{table:rhoEx2}  
 \end{table}

\subsection{Example 3}\label{ex3}
\begin{figure}[!hbt]
\begin{minipage}{0.45\textwidth}
   	\centering
   	\includegraphics[width=1.00\textwidth]{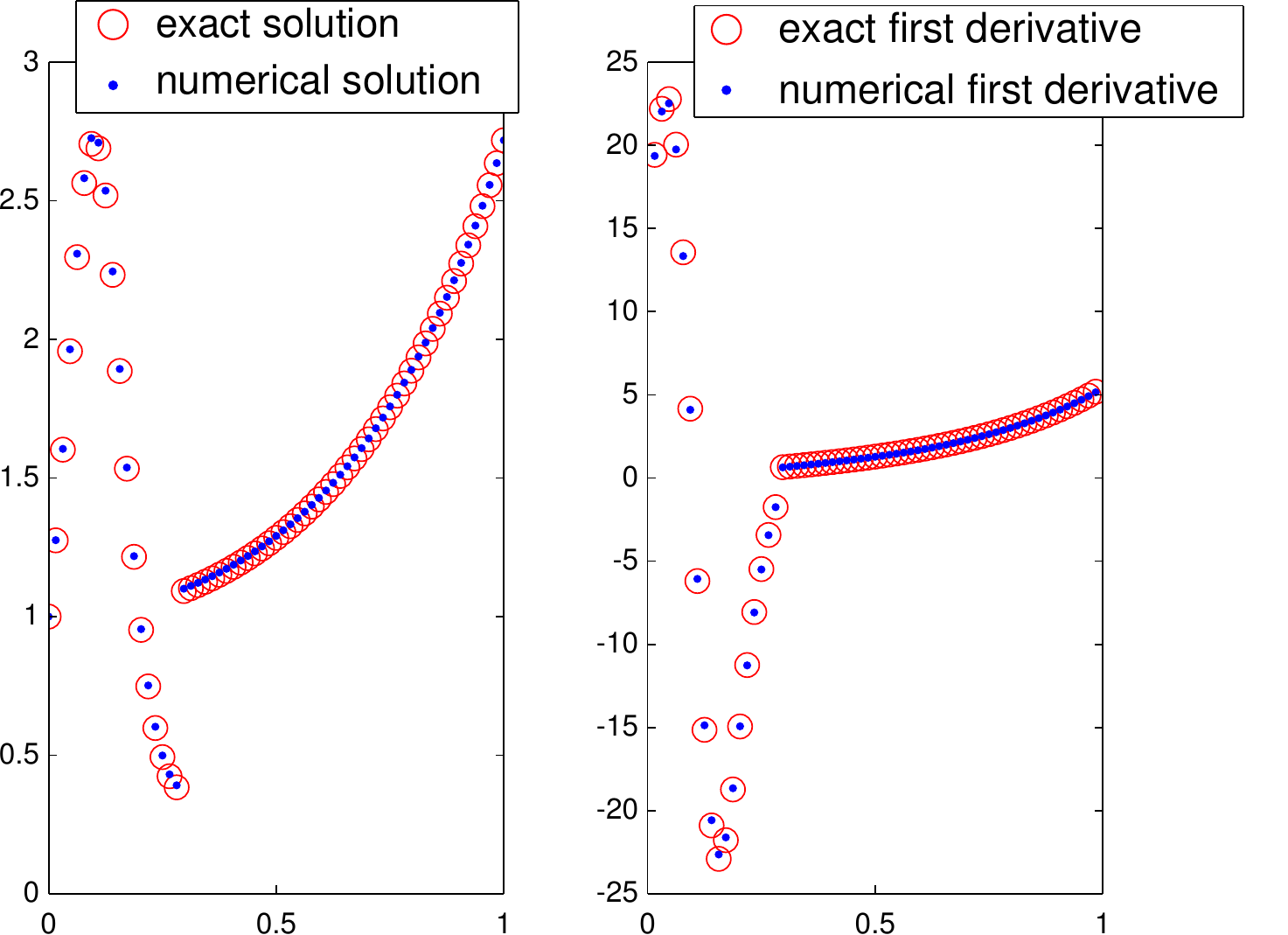}
 \end{minipage}
 \ \hspace{2mm} 
 \begin{minipage}{0.45\textwidth}
  	\centering
  	\captionsetup{width=0.80\textwidth}
	\includegraphics[width=1.00\textwidth]{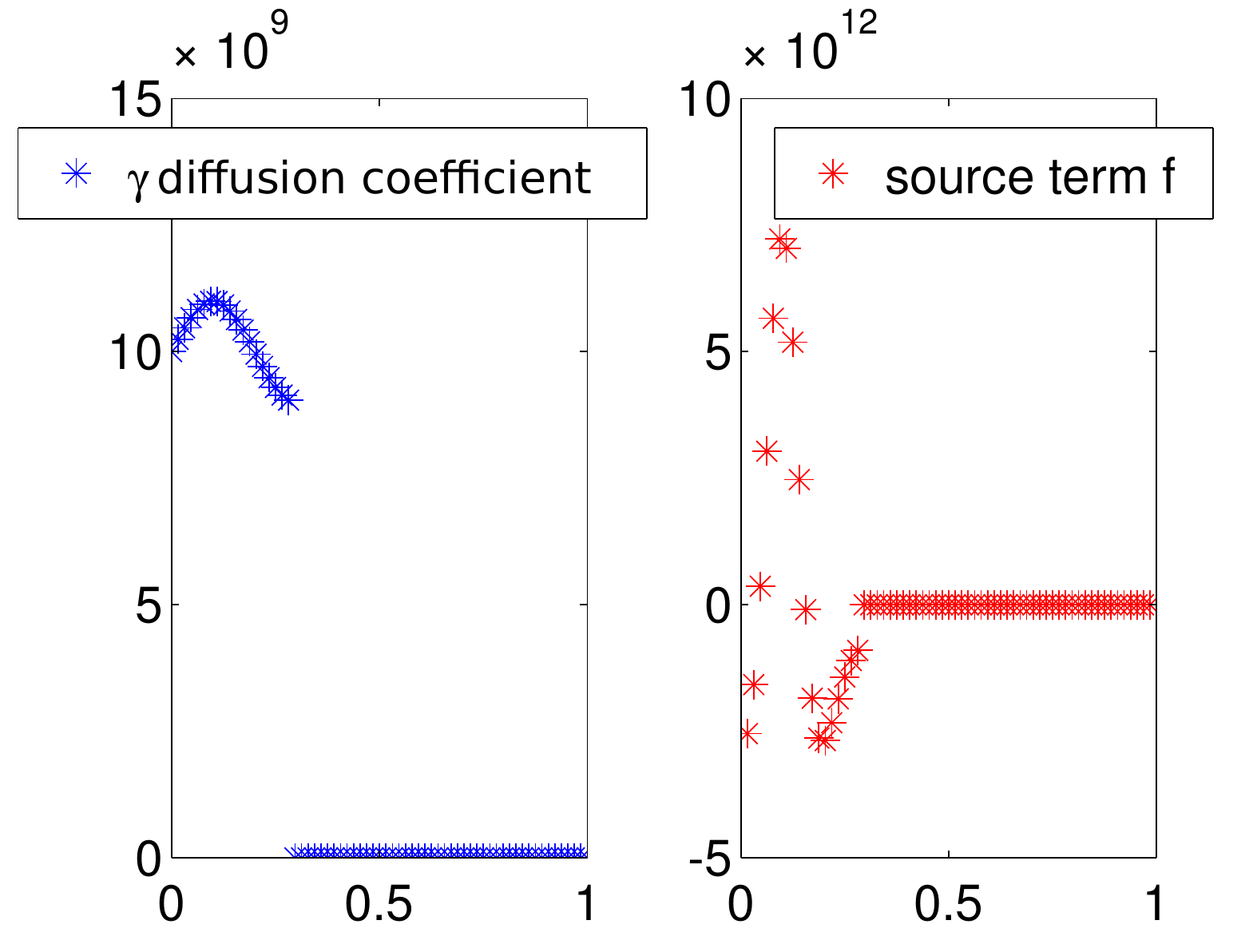}
 \end{minipage}
 \captionsetup{width=0.80\textwidth}
 \caption{\footnotesize{We refer to Ex. \ref{ex3}. The data are computed for $N=64$.}}
 	\label{fig:dataEx3}
\end{figure}

\begin{figure}[!hbt]
 \captionsetup{width=0.80\textwidth}
 \begin{minipage}{0.45\textwidth}
 \centering
   	\includegraphics[width=1.00\textwidth]{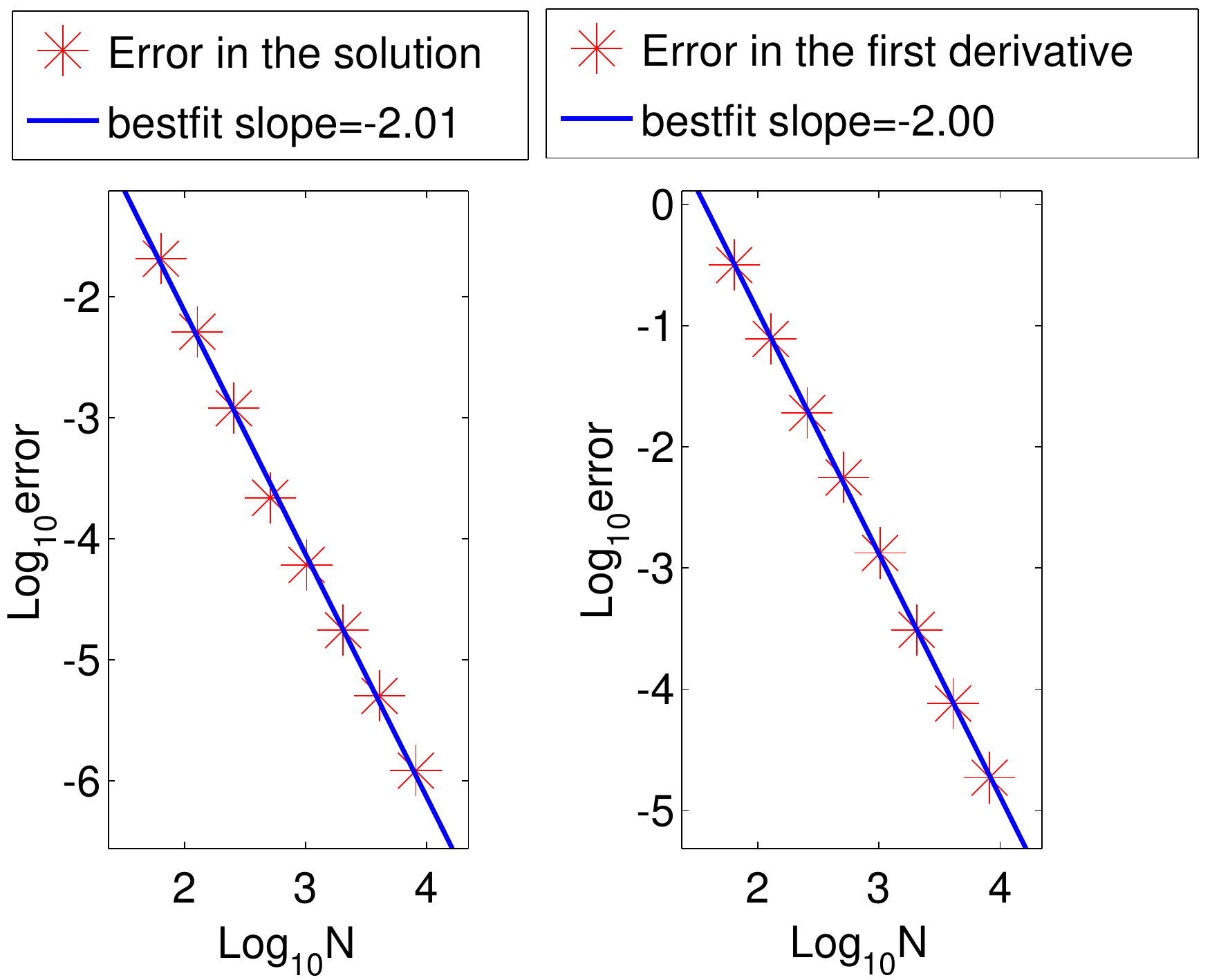}
 \end{minipage}
   \hspace{2mm} 
 \begin{minipage}{0.45\textwidth}
 \centering
\begin{tabular}{||c || c | c || c | c ||}  
\hline\hline                        
$N+1$ & $\left\| \vec{u}-\uh \right\|_\infty$ & order & $\left\| \vec{u}'-\uh' \right\|_\infty$ & order \\ [0.5ex] 
\hline                    
64 & 2.07 $\cdot 10^{-2}$ & - & 3.15 $\cdot 10^{-1}$  & -   \\  
128 & 5.15 $\cdot 10^{-3}$ & 2.01 & 7.76 $\cdot 10^{-2}$  & 2.02   \\  
256 & 1.20 $\cdot 10^{-3}$ & 2.10 & 1.90 $\cdot 10^{-2}$  & 2.03   \\  
512 & 2.19 $\cdot 10^{-4}$ & 2.46 & 5.57 $\cdot 10^{-3}$  & 1.77   \\  
1024 & 6.10 $\cdot 10^{-5}$ & 1.84 & 1.33 $\cdot 10^{-3}$  & 2.07   \\  
2048 & 1.76 $\cdot 10^{-5}$ & 1.79 & 3.09 $\cdot 10^{-4}$  & 2.11   \\  
4096 & 5.05 $\cdot 10^{-6}$ & 1.80 & 7.60 $\cdot 10^{-5}$  & 2.02   \\  
8192 & 1.22 $\cdot 10^{-6}$ & 2.05 & 1.86 $\cdot 10^{-5}$  & 2.03   \\  
\hline     
 \end{tabular} 
 \end{minipage}
	\caption{\footnotesize{We refer to Ex \ref{ex3}.
	Left: Representation of the $L^\infty$-error of the numerical solution and its derivative. The slope of the best-fit lines is respectively $s=-2.01$ and $s=-2.00$.
	 Right: List of errors and order of accuracy computed by subsequent errors. 
	}}
	\label{fig:bestfitEx3}
 	\end{figure}

We choose (see Fig. \ref{fig:dataEx3})
\[
\alpha=0.283 \; \; \;
\left\{
\begin{array}{ccc}
u^L&=&e^{\sin(5 \pi x)} \\
u^R&=&e^{x^2}
\end{array},
\right. \; \; \;
\left\{
\begin{array}{ccc}
\gamma^L&=&10^9\left( 10+\sin(5 \pi x) \right) \\
\gamma^R&=&3+\cos(5 \pi x)
\end{array}.
\right.
\]
Fig. \ref{fig:bestfitEx3} shows the numerical results and the second order slope of the best-fit line for the $L^\infty$-error of the numerical solution and its derivative. 
Table \ref{table:rhoEx3} shows the convergence factor for different values of $N$ and $N_c$.

\begin{table}[!hbt]
\captionsetup{width=0.80\textwidth}
\caption{ \footnotesize{ Measured $V(1,1)$-cycle convergence factor for the numerical test of Ex. \ref{ex3}. We use $N+2$ number of grid points in the finest grid; $N_c+2$ number of grid points in the coarsest grid.  }} 
\centering      
\begin{tabular}{|| c c || c | c | c | c | c | c | c | c ||}  
\hline\hline                        
 & N+1 & 32 & 64 & 128 & 256 & 512 & 1024 & 2048 & 4096  \\ [0.5ex] 
Nc+1 & &  &  &  &  &  &  &  & \\ 
\hline\hline                        
16 & & 0.09 & 0.10 & 0.12 & 0.15 & 0.15 & 0.15 & 0.15 & 0.15\\ 
32 & &  & 0.09 & 0.10 & 0.15 & 0.15 & 0.15 & 0.15 & 0.15\\ 
64 & &  &  & 0.12 & 0.15 & 0.15 & 0.15 & 0.15 & 0.15\\ 
128 & &  &  &  & 0.13 & 0.15 & 0.15 & 0.15 & 0.15\\ 
\hline\hline                        
 \end{tabular} 
      \label{table:rhoEx3}  
 \end{table}

\subsection{Example 4}\label{ex4}
\begin{figure}[!hbt]
\begin{minipage}{0.45\textwidth}
   	\centering
   	\includegraphics[width=1.00\textwidth]{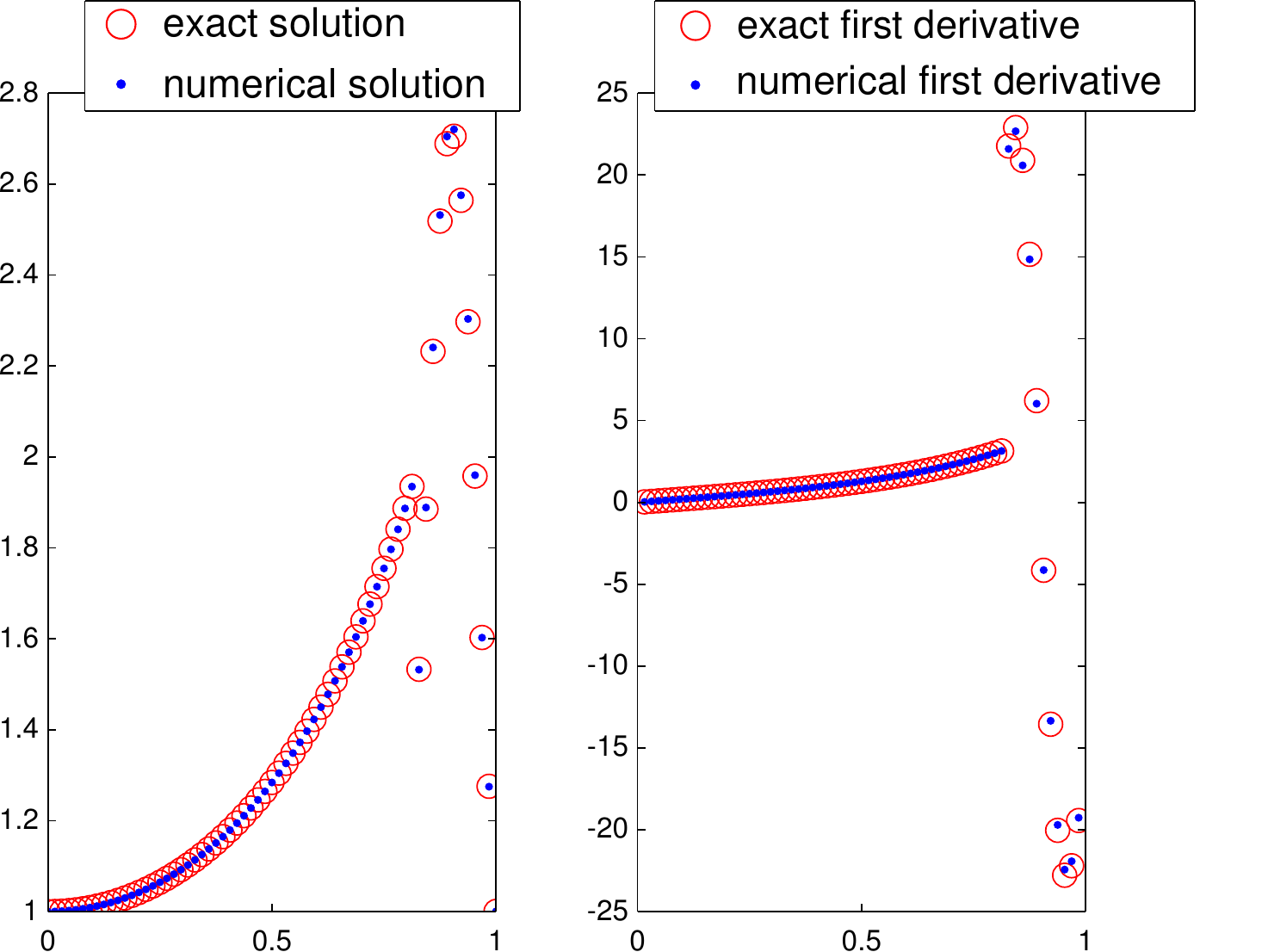}
 \end{minipage}
 \ \hspace{2mm} 
 \begin{minipage}{0.45\textwidth}
  	\centering
  	\captionsetup{width=0.80\textwidth}
	\includegraphics[width=1.00\textwidth]{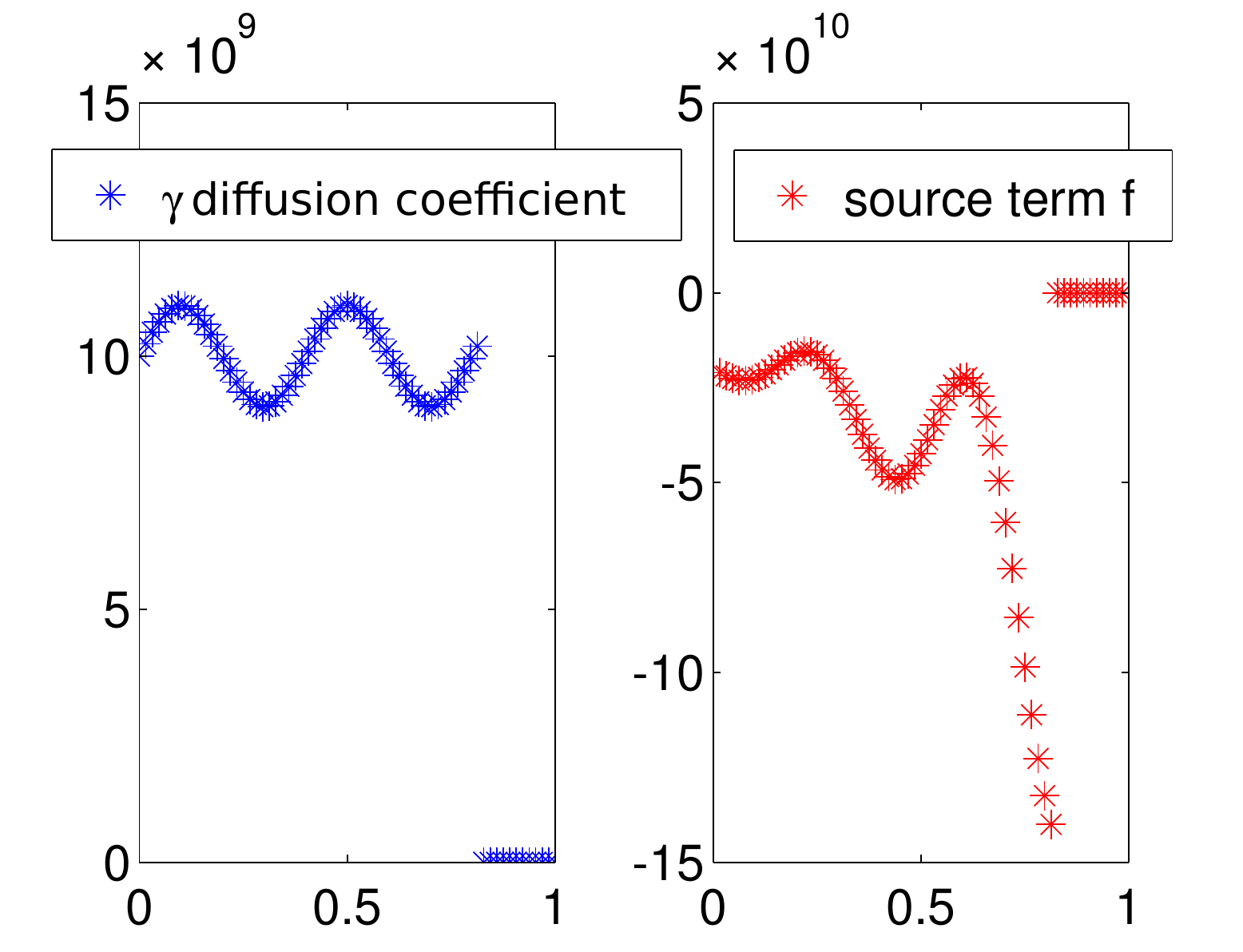}
 \end{minipage}
 \captionsetup{width=0.80\textwidth}
 \caption{\footnotesize{We refer to Ex. \ref{ex4}. The data are computed for $N=64$.}}
 	\label{fig:dataEx4}
\end{figure}

\begin{figure}[!hbt]
 \captionsetup{width=0.80\textwidth}
 \begin{minipage}{0.45\textwidth}
 \centering
   	\includegraphics[width=1.00\textwidth]{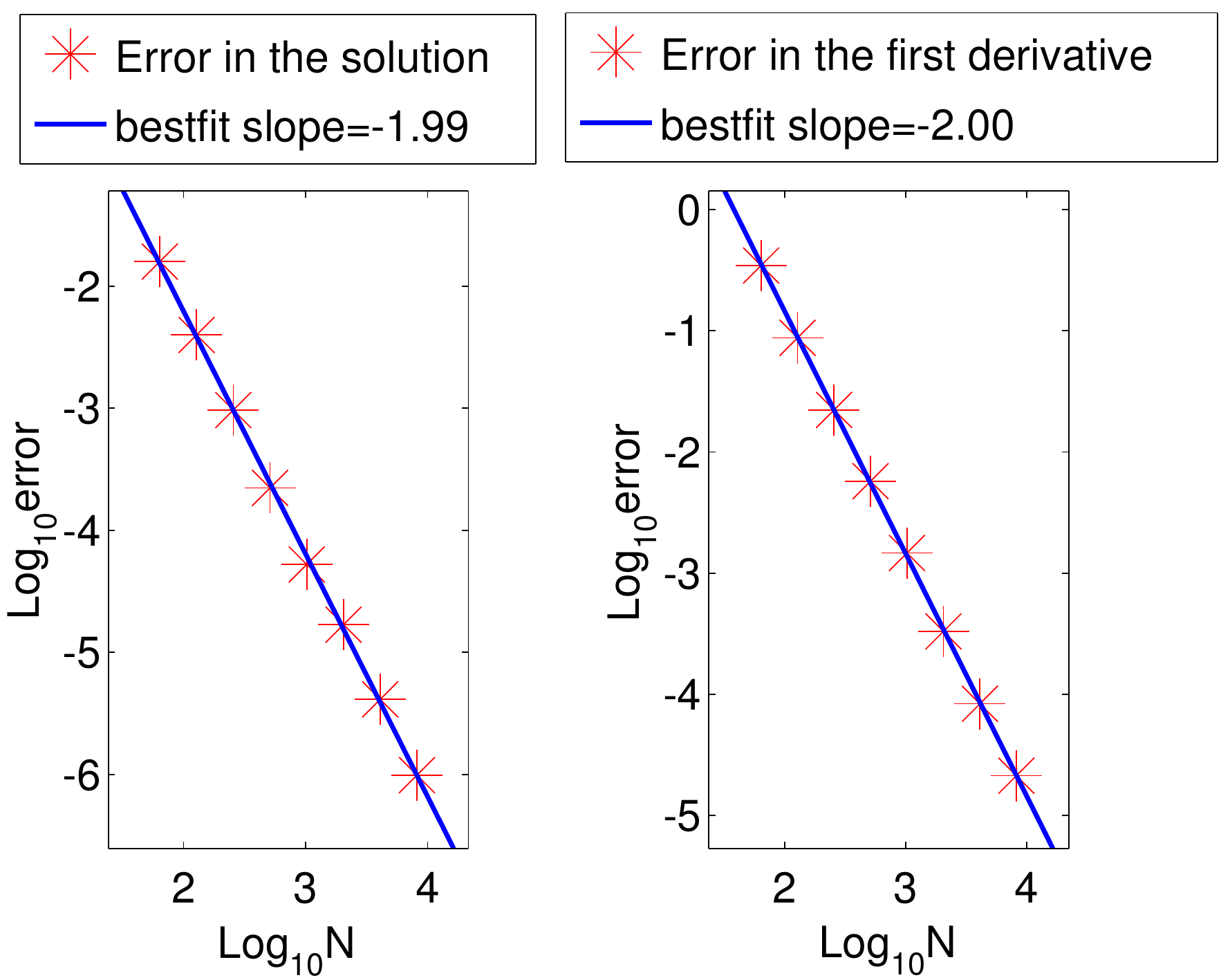}
 \end{minipage}
   \hspace{2mm} 
 \begin{minipage}{0.45\textwidth}
 \centering
\begin{tabular}{||c || c | c || c | c ||}  
\hline\hline                        
$N+1$ & $\left\| \vec{u}-\uh \right\|_\infty$ & order & $\left\| \vec{u}'-\uh' \right\|_\infty$ & order \\ [0.5ex] 
\hline                    
64 & 1.59 $\cdot 10^{-2}$ & - & 3.46 $\cdot 10^{-1}$  & -   \\  
128 & 3.99 $\cdot 10^{-3}$ & 2.00 & 8.74 $\cdot 10^{-2}$  & 1.98   \\  
256 & 9.66 $\cdot 10^{-4}$ & 2.05 & 2.22 $\cdot 10^{-2}$  & 1.98   \\  
512 & 2.23 $\cdot 10^{-4}$ & 2.12 & 5.74 $\cdot 10^{-3}$  & 1.95   \\  
1024 & 5.25 $\cdot 10^{-5}$ & 2.08 & 1.47 $\cdot 10^{-3}$  & 1.97   \\  
2048 & 1.68 $\cdot 10^{-5}$ & 1.64 & 3.31 $\cdot 10^{-4}$  & 2.15   \\  
4096 & 4.12 $\cdot 10^{-6}$ & 2.03 & 8.36 $\cdot 10^{-5}$  & 1.98   \\  
8192 & 9.87 $\cdot 10^{-7}$ & 2.06 & 2.13 $\cdot 10^{-5}$  & 1.97   \\  
\hline     
 \end{tabular} 
 \end{minipage}
	\caption{\footnotesize{We refer to Ex \ref{ex4}.
	Left: Representation of the $L^\infty$-error of the numerical solution and its derivative. The slope of the best-fit lines is respectively $s=-1.99$ and $s=-2.00$.
	 Right: List of errors and order of accuracy computed by subsequent errors. 
	}}
	\label{fig:bestfitEx4}
 	\end{figure}

We choose (see Fig. \ref{fig:dataEx4})
\[
\alpha=0.813, \; \; \;
\left\{
\begin{array}{ccc}
u^L&=&e^{x^2} \\
u^R&=&e^{\sin(5 \pi x)}
\end{array},
\right. \; \; \;
\left\{
\begin{array}{ccc}
\gamma^L&=&10^9\left( 10+\sin(5 \pi x) \right) \\
\gamma^R&=&3+\cos(5 \pi x)
\end{array}.
\right.
\]
Fig. \ref{fig:bestfitEx4} shows the numerical results and the second order slope of the best-fit line for the $L^\infty$-error of the numerical solution and its derivative. 
Table \ref{table:rhoEx4} shows the convergence factor for different values of $N$ and $N_c$.

\begin{table}[!hbt]
\captionsetup{width=0.80\textwidth}
\caption{ \footnotesize{ Measured $V(1,1)$-cycle convergence factor for the numerical test of Ex. \ref{ex4}. We use $N+2$ number of grid points in the finest grid; $N_c+2$ number of grid points in the coarsest grid.  }} 
\centering      
\begin{tabular}{|| c c || c | c | c | c | c | c | c | c ||}  
\hline\hline                        
 & N+1 & 32 & 64 & 128 & 256 & 512 & 1024 & 2048 & 4096  \\ [0.5ex] 
Nc+1 & &  &  &  &  &  &  &  & \\ 
\hline\hline                        
16 & & 0.17 & 0.12 & 0.14 & 0.18 & 0.17 & 0.15 & 0.16 & 0.15\\ 
32 & &  & 0.11 & 0.14 & 0.16 & 0.15 & 0.15 & 0.15 & 0.15\\ 
64 & &  &  & 0.06 & 0.14 & 0.15 & 0.15 & 0.15 & 0.15\\ 
128 & &  &  &  & 0.12 & 0.15 & 0.15 & 0.15 & 0.15\\ 
\hline\hline                        
 \end{tabular} 
      \label{table:rhoEx4}  
 \end{table}

\subsection{Independence of convergence factor from the jump in the coefficient}
In this section we show that the convergence factor does not depend on the jump in the coefficient.
We choose
\[
\alpha=0.543, \; \; \;
\left\{
\begin{array}{ccc}
u^L&=&0 \\
u^R&=&0
\end{array},
\right. \; \; \;
\left\{
\begin{array}{ccc}
\gamma^L&=&10^p \\
\gamma^R&=&1
\end{array}
\right.
\]
and start the multigrid process with an initial guess different from zero, in order to compute the asymptotic convergence factor. We list the results in Table \ref{table:CVindep}.

\begin{table}[!hbt]
\captionsetup{width=0.80\textwidth}
\caption{ \footnotesize{ Measured $V(1,1)$ asymptotic convergence factors for a problem with a jumping coefficient of the order $10^p$  }} 
\centering      
\begin{tabular}{|| l || c | c | c | c | c | c ||}  
\hline\hline                        
 $p$ & 0 & 1 & 2 & 3 & 4 & 5  \\ [0.5ex] 
\hline\hline                        
$\rho$ & 0.11 & 0.10 & 0.11 & 0.11 & 0.11 & 0.10 \\
\hline\hline                        
 \end{tabular} 
      \label{table:CVindep}  
 \end{table}  
 
\textbf{Remark. (Comparison with Domain Decomposition Method)} 
\\Domain Decomposition Method (DDM) is another iterative method to solve elliptic problems with discontinuous coefficient, based on solving iteratively the two subproblems
\begin{equation}\label{DDsub1}
\left\{
\begin{array}{rcll}
- \pad{}{x} \left( \gamma^L \pad{u^{L,(m+1)}}{x} \right) &=& f & \mbox{ in } [0,\alpha[ \\
u^{L,(m+1)}(0)&=&g_0 & \\
u^{L,(m+1)}(\alpha)&=&u^{R,(m)}(\alpha) &
\end{array}
\right.
\end{equation}
\begin{equation}\label{DDsub2}
\left\{
\begin{array}{rcll}
- \pad{}{x} \left( \gamma^R \pad{u^{R,(m+1)}}{x} \right) &=& f & \mbox{ in } ]\alpha,1] \\
\gamma^R \: \pad{u^{R,(m+1)}(\alpha)}{x}&=&\gamma^L \: \pad{u^{L,(m+1)}(\alpha)}{x} & \\
u^{R,(m+1)}(1)&=&g_1 &
\end{array}
\right.
\end{equation}
until convergence.
A little drawback of this method is that, in order to guarantee the convergence, it must be $\alpha>0.5$ (see~\cite[pag. 12]{QuarteroniValli:DDM}). Our method may be regarded as a DDM, but in place of solving a subproblem to provide the right-hand side for the other subproblem (and so on iteratively), we just perform a relaxation on a subproblem, and with the guess obtained we build the right-hand side of the other subproblem, as it can be seen in Sec. \ref{ItMethod}. With this relaxing strategy, the convergence is always guaranteed, as showed in numerical tests.

\section*{Conclusion}
A second order discretization for elliptic equation with discontinuous coefficient on an arbitrary interface has been provided. Second order accuracy in the derivative is obtained as well.
The linear system is solved by an iterative method obtained relaxing the interface conditions.
The iterative method is then speeded up by a proper multigrid approach, which transfers separately the defect for both sub-problems obtained from the multi-domain formulation. The measured convergence factor is close to the one measured in the case of smooth coefficients and it does not depend on the magnitude of the jump in the coefficient.
The method is similar to Domain Decomposition Methods, but a single relaxation sweep is performed in each subdomain instead to solve it completely. This makes the method more flexible and there is no restriction on the relative size of the two subdomains.
\\This paper is the building-block for a future work in higher dimension~\cite{CocoRusso:discontinuous2d}, which will be carried out by combining the second order discretization in arbitrary domain with smooth coefficients~\cite{CocoRusso:Elliptic} and the multigrid treatment of problems with non-eliminated boundary conditions in arbitrary domain~\cite{CocoRusso:MG}.
\\Other future works concern the convection-diffusion equation in a moving domain, in order to study applications modeled by a Stefan-Type problem. A level-set function will keep track of the moving interface.
\\All this extensions will be coupled with the use of Adaptive Mesh Refinement to obtain accurate solution in the case of domain with complex boundary. A proper multigrid approach is under investigation for all these works.

\addcontentsline{toc}{section}{References}
\bibliographystyle{abbrv}
\bibliography{bibliography}

\end{document}